\documentclass{svjour3}                     
\smartqed  
\usepackage{graphicx}
\usepackage{amssymb}

\usepackage{amsfonts,amsmath,latexsym}
\usepackage[active]{srcltx}
\usepackage{color}

\usepackage{enumerate}
\usepackage{stackrel}
\usepackage{stackengine}
\usepackage{mathtools}
\usepackage{verbatim}

\DeclareMathOperator{\Span}{span}

\newcommand{\bu}{\mbox{\boldmath{$u$}}}
\newcommand{\bt}{\mbox{\boldmath{$t$}}}
\newcommand{\bcero}{\mbox{\boldmath{$0$}}}

\newcommand{\bx}{\mbox{\boldmath{$\bx$}}}
\newcommand{\by}{\mbox{\boldmath{$\by$}}}
\newcommand{\bv}{\mbox{\boldmath{$v$}}}
\newcommand{\bw}{\mbox{\boldmath{$w$}}}
\newcommand{\bz}{\mbox{\boldmath{$v$}}}

\newcommand{\var}{\varepsilon}

\newcommand{\ba}{\mbox{\boldmath{$a$}}}

\newcommand{\fb}{\mbox{\boldmath{$f$}}}
\newcommand{\bg}{\mbox{\boldmath{$g$}}}
\renewcommand{\by}{\mbox{\boldmath{$y$}}}
\renewcommand{\bx}{\mbox{\boldmath{$x$}}}
\newcommand{\be}{\mbox{\boldmath{$e$}}}
\newcommand{\bn}{\mbox{\boldmath{$n$}}}
\newcommand{\bh}{\mbox{\boldmath{$h$}}}
\newcommand{\bsigma}{\mbox{\boldmath{$\sigma$}}}
\newcommand{\btheta}{\mbox{\boldmath{$\theta$}}}

\newcommand{\beeta}{\mbox{\boldmath{$\eta$}}}
\newcommand{\bxi}{\mbox{\boldmath{$\xi$}}}

\newcommand{\bzero}{\mbox{\boldmath{$0$}}}
\newcommand{\bTheta}{\mbox{\boldmath{$\Theta$}}}
\newcommand{\bUcal}{\mbox{\boldmath{$\mathcal{U}$}}}

\newcommand{\Real}{\mathbb R}

\newcommand{\Gamae}{\Gamma^{\varepsilon }}

\newcommand{\dem}{ \textbf{Proof.}\,\,}
\newcommand{\cqd}
 {\hfill $\sqcup\!\!\!\!\sqcap\bigskip$}

\renewcommand{\d}{\partial}
\newcommand{\eij}{e_{i||j}}
\newcommand{\deij}{\dot{e}_{i||j}}
\newcommand{\ekl}{e_{k||l}}

\newcommand{\eab}{e_{\alpha||\beta}}
\newcommand{\est}{e_{\sigma||\tau}}

\newcommand{\eatres}{e_{\alpha||3}}

\newcommand{\edtres}{e_{3||3}}
\newcommand{\deab}{\dot{e}_{\alpha||\beta}}

\newcommand{\dedtres}{\dot{e}_{3||3}}
\newcommand{\gab}{\gamma_{\alpha\beta}}
\newcommand{\gst}{\gamma_{\sigma\tau}}

\renewcommand{\a}{a^{\alpha\beta\sigma\tau}}

\newcommand{\aeps}{a^{\alpha\beta\sigma\tau,\varepsilon}}

\newcommand{\aes}{\ a.e. \ \textrm{in} \ (0,T)}
\newcommand{\aesth}{\ a.e. \ \textrm{in} \ (0,T-h)}
\newcommand{\forallt}{\ \forall  \ t\in[0,T]}

\newcommand{\deb}{\rightharpoonup}
\newcommand{\en}{ \ \textrm{in} \ }
\newcommand{\on}{ \ \textrm{on} \ }
\newcommand{\into}{\int_{\omega}}
\newcommand{\intO}{\int_{\Omega}}

\newcommand{\ten}{(a^{\alpha \sigma}a^{\beta \tau} + a^{\alpha \tau}a^{\beta\sigma})}

\newcommand{\Vcart}{V(\hat{\Omega}^\var)}
\newcommand{\Scart}{S(\hat{\Omega}^\var)}
\newcommand{\vth}{\vartheta}

\begin{document}

\title{Asymptotic analysis of a problem for dynamic thermoelastic shells
}

\titlerunning{Thermoelastic Elliptic Shells }        

\author{M.T. Cao-Rial \and G.~Casti\~neira \and \'Angel Rodr\'{\i}guez-Ar\'os%
\and  S. Roscani  
}


\institute{M.T. Cao-Rial \at               E.T.S. N\'autica e M\'aquinas,
Paseo de Ronda, 51, 15011, Departamento de M\'etodos Matem\'aticos
e Representaci\'on, Universidade da Coru\~na, Spain \\
              Tel.: +34-981167000\\
              \email{teresa.cao@udc.es}           %
\and G. Casti\~neira \at
Área de Matemáticas, Centro Universitario de la Defensa, Universidade de Vigo,
Escuela Naval Militar,  Plaza de España, s/n 36920 Marín (Pontevedra), Spain\\
\email{gonzalo.castineira.v@cud.uvigo.es}
\and \'A. Rodr\'{\i}guez-Ar\'os \at
              E.T.S. N\'autica e M\'aquinas,
Paseo de Ronda, 51, 15011, Departamento de M\'etodos Matem\'aticos
e Representaci\'on, Universidade da Coru\~na, Spain \\
              \email{angel.aros@udc.es}           
\and S. Roscani \at CONICET - Departamento de. Matemática, FCE, Universidad Austral, Paraguay 1950, S2000FZF, Rosario, Argentina\\
\email{sroscani@austral.edu.ar}
}

\date{Received: date / Accepted: date}

\maketitle

\begin{abstract}

In this paper we consider a family of three-dimensional problems in thermoelasticity for elliptic membrane shells and study the asymptotic behaviour of the solution when the thickness tends to zero. We fully characterize with strong convergence results the limit as the unique solution of a two-dimensional problem, where the reference domain is the common middle surface of the family of three-dimensional shells. The problems are dynamic and the constitutive thermoelastic law is given by the Duhamel-Neumann relation.

\end{abstract}

\section{Introduction}
\label{intro}


In the last decades, asymptotic methods have been used to derive and justify simplified models for three-dimensional solid mechanics problems for beams, plates and shells. The foundation for these methods was established by Lions in \cite{lions} and some of the first applications were to plate bending models in \cite{CD,Destuynder}. Many other models for plates have been justified by using asymptotic methods and a comprehensive review concerning plate models can be found in \cite{Ciarlet3}.

Since then, its application has been extended over the years to many other problems like beam bending, rod stretching and elastic shells. For example, Bernoulli-Navier model was justified in \cite{BermudezViano} and the Saint-Venant, Timoshenko and Vlassov models of elastic beams were justified in \cite{TraViano}, while a model for Kelvin-Voigt viscoelastic beams was justified in \cite{AV2} and models for piezoelectric beams were obtained in \cite{VFRR}, all of them by using the asymptotic expansion method followed  by rigorous convergence results. The asymptotic modelling of rods in linearized thermoelasticity was also studied in \cite{TraViano}.


Regarding elastic shells, a complete theory can be found in \cite{Ciarlet4b}, where models  for elliptic membranes, generalized membranes and flexural shells are presented. In there, the reader can find a full description of the asymptotic procedure that leads to the corresponding sets of two-dimensional equations. More recently, in a series of papers we studied the asymptotic analysis of viscoelastic shells \cite{intro2,eliptico,flexural,koiter} and contact problems for elastic shells \cite{ContactShell,AA_contact_shell,Cao_Aros,ArosCao19}. For the dynamic case, the authors in \cite{Limin_mem,Limin_flex} use the asymptotic analysis to derive two-dimensional sets of equations for elastic membranes and flexural shells, though no strong convergence results are provided. Dynamic problems for shells is a topic which is attracting a considerable effort in modeling, analysis and numerical approximation, due to the abundance of real world applications, see for example \cite{piersanti} and references therein.

The aim of the present paper is to provide the first results of the asymptotic analysis devoted to thermoelastic shells in a dynamic regime. Here we briefly describe the formal asymptotic analysis and the limit two-dimensional problem and we focus in the case of elliptic membrane shells, for which we provide a rigorous convergence result. We also discuss the existence and uniqueness of solution for both the three-dimensional problem and the corresponding two-dimensional limit problem.

The structure of the paper is the following: in Section \ref{problema} we shall describe the variational and mechanical formulations of the problem in cartesian coordinates in a general domain, and present a result of existence and uniqueness of solution for that problem. In Section \ref{seccion_dominio_ind} we consider the particular case when the deformable body is, in fact, a shell and reformulate the variational formulation in curvilinear coordinates. Then we give the scaled formulation. To do that, we will use a projection map into a reference domain and we will introduce the scaled unknowns and forces as well as the assumptions on coefficients. We also devote this section to recall and derive results that will be needed later. In Section \ref{procedure} we briefly describe the formal asymptotic analysis which leads to the formulation of limit two-dimensional problems.  Then, in Section \ref{convergence} we discuss the existence and uniqueness of solution for the two-dimensional limit problem and then we focus on the elliptic membrane case, for which we provide a rigorous convergence result. Finally, in Section \ref{descaling} we show that the solution to the re-scaled version of this problem, with true physical meaning, also converges. The paper ends with Section \ref{conclusiones}, devoted to the conclusions and future work.


\section{A three-dimensional dynamic problem for thermoelastic bodies} \label{problema}

Let $\hat{\Omega}^\var$ be a three-dimensional bounded domain and assume that $\bar{\hat{\Omega}}^\varepsilon$ is the reference configurarion of a deformable body made of an elastic material, which is homogeneous and isotropic, with Lam\'e coefficients $\hat{\lambda}^\var\geq0, \hat{\mu}^\var>0$. Let $\hat{\Gamma}^\var=\partial\hat{\Omega}^\var$ denote the boundary of the body, which is divided into two disjoint parts $\hat{\Gamma}^\varepsilon_N$ and $\hat{\Gamma}_0^\varepsilon$, where the measure of the latter is strictly positive. Let $\hat{\bx}^\varepsilon=(\hat{x}_i^\varepsilon)$ be a generic point of $\bar{\hat{\Omega}}^\varepsilon$.
Notice that at first glance, the notation for sets, variables and functions seems unnecessarily complex. Indeed, the ${}^\varepsilon$ and $\hat{ }$ marks are only meaningful in the context of the shells setting, to be detailed in the forthcoming sections. But, given that there we are going to recall results and arguments developed in this current section, we decided to keep here this notation, in favor of future coherence.

We suppose that the material has a thermal dilatation coefficient $\hat{\alpha}_T^\var$, a thermal conductivity coefficient $\hat{k}^\var$, a specific heat coefficient $\hat{\beta}^\var$ and a specific mass density $\hat{\rho}^\var$.  The constitutive equation relating the stress tensor components $\hat{\sigma}_{ij}^\var$ to the linearized strain tensor $\hat{e}_{ij}^\var(\hat{\bu}^\var)$ components, and the temperature $\hat{\vartheta}^\var$ is given by the linearized Duhamel-Neumann law (see, for example \cite{TraViano} and references therein):
\begin{equation}\label{Duhamel-Neumann}
  \hat{\sigma}_{ij}^\var=\hat{\lambda}^\var \hat{e}_{kk}^\var(\hat{\bu}^\var)\delta_{ij}+2\hat{\mu}^\var \hat{e}_{ij}^\var(\hat{\bu}^\var)-\hat{\alpha}_T^\var(3\hat{\lambda}^\var+2\hat{\mu}^\var)\hat{\vartheta}^\var\delta_{ij},
\end{equation}
where $\hat{e}^\var_{ij}(\hat{\bv}^\var)=\frac{1}{2}(\hat{\partial}_j\hat{v}_i^\var+\hat{\partial}_i\hat{v}_j^\var)$ denotes the deformation operator. Here $\delta_{ij}$ represents the Kronecker's symbol and $\hat{\partial}_i$ the partial derivative with respect to $\hat{x}_i^\var$.
Notice that here and below, and for the sake of a clearer exposition, we shall omit the explicit dependence of the various functions on space and time variables, as long as there is no ambiguity.
We assume that the body is subjected to a boundary condition of place; in particular, the displacements field vanishes on  $\hat{\Gamma}_0^\varepsilon$.

The body is under the effect of a heat source $\hat{q}^\var$ and applied volume forces of density $\hat{\fb}^\var=(\hat{f}^{i,\var})$ acting in $\hat{\Omega}^\var$, and tractions of density $\hat{\bh}^\var=(\hat{h}^{i,\var})$ acting upon $\hat{\Gamma}^\var_N$.

Then, the set of equations describing the mechanical behaviour of a regular three-dimensional deformable solid in thermoelasticity are the following:
\begin{problem}\label{problema_fuerte}
Find the displacements field $\hat{\bu}^\varepsilon=(\hat{u}^\varepsilon_i)$ and the temperature field $\hat{\vartheta}^\varepsilon$ verifying
\begin{align*}
&  \hat{\rho}^\var \ddot{\hat{\bu}}^\var-{\rm div}\hat{\bsigma}^\var=\hat{\fb}^\var & \mbox{ in } \hat{\Omega}^\var\times(0,T) , \\
&  \hat{\beta}^\var \dot{\hat{\vartheta}}^\varepsilon=\partial_j(\hat{k}^\varepsilon\hat{\partial}_j^\varepsilon\hat{\vartheta}^\varepsilon)-
\hat{\alpha}_T^\var(3\hat{\lambda}^\var+2\hat{\mu}^\var)\hat{e}^\varepsilon_{kk}(\dot{\hat{\bu}}^\var)+\hat{q}^\var &  \mbox{ in } \hat{\Omega}^\var\times(0,T),\\%
&  \hat{\bu}^\var=\bzero &  \mbox{ on } \hat{\Gamma}_0^\var\times(0,T), \\
&  \hat{\vartheta}^\varepsilon=0 & \mbox{ on } \hat{\Gamma}_0^\var\times(0,T),  \\
&  \hat{\sigma}^\var\hat{\bn}^\var=\hat{\bh}^\var & \mbox{ on } \hat{\Gamma}_N^\var\times(0,T), \\
& \hat{k}^\varepsilon\hat{\partial}_j^\varepsilon\hat{\vartheta}^\varepsilon n_j=0  & \mbox{ on } \hat{\Gamma}_N^\var
\times(0,T), \\
& \dot{\hat{\bu}}^\varepsilon(\cdot,0)=\hat{\bu}^\varepsilon(\cdot,0)=\bzero & \mbox{ in } \hat{\Omega}^\var, \\
&  \hat{\vartheta}^\var(\cdot,0)=0 &  \mbox{ in } \hat{\Omega}^\var,
\end{align*}
where $\hat{\bsigma}^\var=(\hat{\sigma}_{ij}^\var)$ is described in (\ref{Duhamel-Neumann}).
\end{problem}
\begin{remark}
We only consider homogeneous boundary and initial conditions for the sake of simplicity. Besides, our main interest is in the asymptotic analysis which follow in the sections below, and nonhomogenous initial conditions do not introduce major additional difficulties in that sense. 
\end{remark}
Now, to derive the variational formulation of the problem, we first define the space of admissible displacements
$$
V(\hat{\Omega}^\var):=\{\hat{\bv}^\varepsilon=(\hat{v}_i^\varepsilon)\in [H^1(\hat{\Omega}^\var)]^3; \hat{\bv}^\varepsilon=\mathbf{\bcero} \ {\rm on} \ \hat{\Gamma}_0^\var\},
$$
which is a Hilbert space equipped with the norm
$$
\|\hat{\bv}^\var\|_{V(\hat{\Omega}^\var)}=\int_{\hat{\Omega}^\var}\hat{e}^\varepsilon_{ij}(\hat{\bv}^\varepsilon)
   \hat{e}^\varepsilon_{ij}(\hat{\bv}^\varepsilon) d\hat{x}^\varepsilon.
$$
We also define the space of admissible temperatures
$$
{S}(\hat{\Omega}^\var):=\{\hat{\varphi}^\varepsilon\in H^1(\hat{\Omega}^\var); \hat{\varphi}^\varepsilon=0 \ {\rm on} \ \hat{\Gamma}_0^\var\},
$$
which is a Hilbert space equipped with the norm
$$
\|\hat{\varphi}^\var\|_{S(\hat{\Omega}^\var)}=\int_{\hat{\Omega}^\var}\hat{\partial}_j^\varepsilon\hat{\varphi}^\varepsilon \hat{\partial}_j^\varepsilon\hat{\varphi}^\var d\hat{x}^\varepsilon.
$$
Besides, as long as there is no room for confusion, we shall avoid specifying the domain in the subindices for the norms' notation. Further, for the sake of simplicity in the formulations to come, we define the following operators:
\begin{itemize}
  \item The bilinear, continuous and coercive forms
\begin{align*}
a^{V,\var}\colon & V(\hat{\Omega}^\var) \times V(\hat{\Omega}^\var) \rightarrow \Real \\
& (\hat{\bu}^\var, \hat{\bv}^\var)   \rightarrow a^{V,\var}(\hat{\bu}^\var, \hat{\bv}^\var)=\int_{\hat{\Omega}^\var}\hat{A}^{ijkl,\varepsilon}\hat{e}^\varepsilon_{kl}(\hat{\bu}^\varepsilon)
   \hat{e}^\varepsilon_{ij}(\hat{\bv}^\varepsilon) d\hat{x}^\varepsilon,\\
a^{S,\var}\colon & S(\hat{\Omega}^\var) \times S(\hat{\Omega}^\var) \rightarrow \Real \\
& (\hat{\varphi}^\var, \hat{\psi}^\var)   \rightarrow a^{S,\var}(\hat{\varphi}^\var, \hat{\psi}^\var)=\int_{\hat{\Omega}^\var}\hat{k}^\varepsilon\hat{\partial}_j^\varepsilon\hat{\varphi}^\varepsilon \hat{\partial}_j^\varepsilon\hat{\psi}^\var d\hat{x}^\varepsilon,
\end{align*}
where $\hat{A}^{ijkl,\varepsilon}=\hat{\lambda}^\var\delta^{ij}\delta^{kl}+\hat{\mu}^\var(\delta^{ik}\delta^{jl}+\delta^{il}\delta^{jk})$
denotes the elasticity fourth-order tensor.
  \item The continuous form
\begin{align*}
c^\var\colon &S(\hat{\Omega}^\var) \times V(\hat{\Omega}^\var) \rightarrow \Real \\
& (\hat{\varphi}^\var, \hat{\bv}^\var)   \rightarrow c^\var(\hat{\varphi}^\var, \hat{\bv}^\var)=\int_{\hat{\Omega}^\var}\hat{\alpha}_T^\var(3\hat{\lambda}^\var+2\hat{\mu}^\var)\hat{\varphi}^\var \hat{e}^\varepsilon_{kk}(\hat{\bv}^\varepsilon)  d\hat{x}^\varepsilon.
\end{align*}
%

%
\item The functional $\hat{J}^\var(\cdot)$ is defined a.e. in $(0,T)$ as
\begin{align*}
\left\langle \hat{J}^\var(t), \hat{\bv}^\var \right\rangle
=\int_{\hat{\Omega}^\var}\hat{f}^{i,\var}(t)\hat{v}_i^\varepsilon\,d\hat{x}^\varepsilon%
+\int_{\hat{\Gamma}_N^\var}\hat{h}^{i,\var}(t)\hat{v}_i^\varepsilon\,d\hat{\Gamma}^\var,\, \qquad \forall \, \hat{\bv}^\var \in  V(\hat{\Omega}^\var),
\end{align*}
where we use the notation for a duality pair $\left\langle \cdot,\cdot\right\rangle$ in $V'(\hat{\Omega}^\var) \times V(\hat{\Omega}^\var)$, and similarly, 
\begin{align*}
\left\langle \hat{Q}^\var(t),\hat{\varphi}^\var\right\rangle = \int_{\hat{\Omega}^\var}\hat{q}^{\var}(t)\hat{\varphi}^\var\,d\hat{x}^\varepsilon  \qquad \forall \, \hat{\varphi}^\var \in S(\hat{\Omega}^\var).
\end{align*}
\end{itemize}

Then, it is straightforward to obtain the following variational formulation:

\begin{problem}\label{problema_cart}
Find a pair $t\mapsto(\hat{\bu}^\varepsilon(\hat{\bx}^\var,t),\hat{\vartheta}^\varepsilon(\hat{\bx}^\var,t))$ of $[0,T]\to V(\hat{\Omega}^\var)\times {S}(\hat{\Omega}^\var)$ verifying
\begin{align}
&\hat{\rho}^\var \left\langle\ddot{\hat{u}}_i^\var,\hat{v}_i^\varepsilon \right\rangle%
+ a^{V,\var}(\hat{\bu}^\var, \hat{\bv}^\var)%
-c^\var(\hat{\vartheta}^\var, \hat{\bv}^\var)%
=\left\langle \hat{J}^\var(t), \hat{\bv}^\var \right\rangle%
   \ \forall \hat{\bv}^\varepsilon\in V(\hat{\Omega}^\var),\aes,\label{P2.2a}\\
&
\hat{\beta}^\var\left\langle \dot{\hat{\vartheta}}^\varepsilon, \hat{\varphi}^\varepsilon \right\rangle
+a^{S,\var}(\hat{\vartheta}^\var, \hat{\varphi}^\var)%
+c^\var(\hat{\varphi}^\var, \dot{\hat{\bu}}^\var)
=\left\langle \hat{Q}^\var(t),\hat{\varphi}^\var\right\rangle   \ \forall \hat{\varphi}^\varepsilon\in {S}(\hat{\Omega}^\var),\aes,\label{P2.2b}
\end{align}
with $\dot{\hat{\bu}}^\varepsilon(\cdot,0)=\hat{\bu}^\varepsilon(\cdot,0)=\bzero$ and $\hat{\vartheta}^\var(\cdot,0)=0$.
\end{problem}

In favour of simplicity, we are going to assume that the different parameters of the problem (thermal conductivity, thermal dilatation, specific heat coefficient, mass density, Lamé coefficients) are constants.

\begin{theorem}\label{Th_exist_unic_3D_bid_cero}
Let us assume that
\begin{equation*}
\begin{cases}
\hat{\fb}^\var \in H^1(0,T;[L^2(\hat{\Omega}^\var)]^3),\\
\hat{\bh}^\var \in H^2(0,T;[L^2(\hat{\Gamma}_N^\var)]^3), \,\mbox{and}\, \, \hat{\bh}^\var(\cdot, 0)=\bzero, \\
\hat{q}^\var \in H^1(0,T;L^2(\hat{\Omega}^\var)).
\end{cases}
\end{equation*}
Then, there exists a unique pair $(\hat{\bu}^\var(\bx,t), \hat{\vth}^\var(\hat{\bx},t))$ solution to Problem \ref{problema_cart} such that
\begin{equation}\label{reg_despl}
\begin{cases}
\hat{\bu}^\var \in L^\infty(0,T;V(\hat{\Omega}^\var))\\
\dot{\hat{\bu}}^\var \in L^\infty(0,T;[L^2(\hat{\Omega}^\var)]^3)\cap L^\infty(0,T;V(\hat{\Omega}^\var)),\\
\ddot{\hat{\bu}}^\var \in L^\infty(0,T;V'(\hat{\Omega}^\var))\cap L^\infty(0,T;[L^2(\hat{\Omega}^\var)]^3),
\end{cases}
\end{equation}

\begin{equation}\label{reg_temp}
\begin{cases}
\hat{\vartheta}^\var \in L^\infty(0,T;L^2(\hat{\Omega}^\var))\cap L^2(0,T;S(\hat{\Omega}^\var)),\\
\dot{\hat{\vartheta}}^\var \in L^\infty(0,T;L^2(\hat{\Omega}^\var))\cap L^2(0,T;S(\hat{\Omega}^\var)).
\end{cases}
\end{equation}
\end{theorem}
\begin{remark}
The regularity results in (\ref{reg_despl}c) and (\ref{reg_temp}b) imply that the duality products involving $\ddot{\hat{\bu}}^\var$ and $\dot{\hat{\vartheta}}^\var$ in \eqref{P2.2a} and \eqref{P2.2b} are actually the usual inner products in $L^2(\hat{\Omega}^\var)$.
\end{remark}
The proof can be derived, for example, by following \cite[p. 359]{Marsden} or \cite{Francfort}. 
We provide below an alternative proof by following the Faedo-Galerkin method.

\proof
Let $\left\{ \hat{\bw}_i \right\}_{i=1}^\infty$ and $\left\{ \hat{s}_i \right\}_{i=1}^\infty$ be two sequences of functions such that

\begin{equation}\label{Vm}
\begin{cases}
\hat{\bw}_i \in \Vcart \quad \forall i,\\
\hat{\bw}_1, \hdots, \hat{\bw}_m \quad \text{are orthonormal functions} \, \text{and }  \, V_m=\left\langle \hat{\bw}_1, \hdots, \hat{\bw}_m \right\rangle,\, \forall m  \\
\Vcart =\overline{\displaystyle\bigcup_{m\geq 1} V_m}.
\end{cases}
\end{equation}

\begin{equation}\label{Sm}
\begin{cases}
\hat{s}_i \in \Scart \quad \forall i,\\
\hat{s}_1, \hdots, \hat{s}_m \quad \text{are orthonormal functions} \, \text{ and }  \, S_m=\left\langle \hat{s}_1, \hdots, \hat{s}_m \right\rangle,\, \forall m  \\
\Scart =\overline{\displaystyle\bigcup_{m\geq 1} S_m}.
\end{cases}
\end{equation}
The approximated solutions $(\hat{\bu}^m,\hat{\vth}^m)$ are defined  by the following problem:
\begin{problem}\label{FG-m}Find the functions  $\hat{\bu}^m\colon [0,T]\rightarrow V_m$  and $\hat{\vth}^m\colon [0,T]\rightarrow S_m$ in the form
\begin{align*}
\hat{\bu}^m(\hat{\bx},t)=\displaystyle\sum_{i=1}^{m}u_i^m(t)\hat{\bw}_i(\hat{\bx}),\\
\hat{\vth}^m(\hat{\bx},t)=\displaystyle\sum_{i=1}^{m}\vth_i^m(t)\hat{s}_i(\hat{\bx}),
\end{align*}
such that
\begin{equation}\label{Galerkin-Eq1}
\hat{\rho}^\var \left\langle \ddot{\hat{\bu}}^m,\hat{\bv}^m \right\rangle
+ a^{V,\var}({\hat{\bu}}^m,\hat{\bv}^m)-c^\var(\hat{\vartheta}^m,\hat{\bv}^m)
=\left\langle \hat{J}^\var(t),\hat{\bv}^m \right\rangle, \qquad \forall \, \hat{\bv}^m \in  V_m,
	\end{equation}
\begin{equation}\label{Galerkin-Eq2}
\hat{\beta}^\var\left\langle \dot{\hat{\vth}}^m,\hat{\varphi}^m  \right\rangle
+ a^{S,\var}({\hat{\vth}}^m,\hat{\varphi}^m)
+ c^\var(\hat{\varphi}^m,\dot{\hat{\bu}}^m)
=\left\langle \hat{Q}^\var(t),\hat{\varphi}^m\right\rangle,  \qquad \forall \, \hat{\varphi}^m \in S_m,
\end{equation}
with the initial conditions
\begin{equation}\label{initial_conditions}
  \hat{\bu}^m(0)=\dot{\hat{\bu}}^m(0)=\bzero,\quad \hat{\vth}^m(0)=0.
\end{equation}
\end{problem}
Finding a solution for Problem \ref{FG-m} is equivalent to solving a first order differential equation system
\begin{equation*}
\dot{\bold{Z}}(t)=\bold{F}(t,\bold{Z}), \quad \bold{Z}(0)=\bold{0}.
\end{equation*}
where $\bold{Z}(t)=({v}^m_1(t),\hdots,{v}_m^m(t),{u}_1^m(t),\hdots,{u}_m^m(t),{\vth}_1^m(t),\hdots,{\vth}_m^m(t))$,
with ${v}_j^m(t)=\dot{{u}}_j^m(t)$.
The Picard-Lindeloff theorem gives a unique absolutely continuous solution in an interval $[0,t_m]$ which depends on the supreme of function $\bold{F}$ (which does not depend on time). Then, being the functions $F_j$  uniformly Lipschitz in the variable $\bold{Z}$, if we prove that the solution $\bold{Z}(t)$ is bounded, we can extend the solution to the whole interval $[0,T]$.\\
\medskip
Now the goal is to obtain estimations in appropriate normed spaces for $\hat{\bu}^m$, $\dot{\hat{\bu}}^m$, $\hat{\vth}^m$ and
$\dot{\hat{\vth}}^m$.

We can take $\hat{\bv}^m=\dot{\hat{\bu}}^m\in V_m$ and $\hat{\varphi}^m=\hat{\vth}^m\in S_m$ in
\eqref{Galerkin-Eq1}, \eqref{Galerkin-Eq2} respectively, and adding both equations we have that

\begin{equation*}
\begin{split}
\hat{\rho}^\var \left\langle \ddot{\hat{\bu}}^m,\dot{\hat{\bu}}^m \right\rangle%
+ a^{V,\var}({\hat{\bu}}^m,\dot{\hat{\bu}}^m)%
+\hat{\beta}^\var\left\langle \dot{\hat{\vth}}^m,\hat{\vth}^m\right\rangle%
+ a^{S,\var}({\hat{\vth}}^m,\hat{\vth}^m)%
\\
=\left\langle \hat{J}^\var(t),\dot{\hat{\bu}}^m \right\rangle + \left\langle \hat{Q}^\var(t),\hat{\vth}^m\right\rangle,
\end{split}
\end{equation*}

or equivalently

\begin{equation}\label{Glaerkin-Eq4}
\begin{split}
\frac{1}{2} \frac{d}{dt}\left\{ \hat{\rho}^\var \left| \dot{\hat{\bu}}^m (t)\right|^2_{0}%
+a^{V,\var}({\hat{\bu}}^m(t), {\hat{\bu}}^m(t))%
+\hat{\beta}^\var\left| \hat{\vth}^m (t)\right|^2_{0} \right\}%
+ a^{S,\var}\left({\hat{\vth}}^m,{\hat{\vth}}^m\right)%
\\%
=\left\langle \hat{J}^\var(t),\dot{\hat{\bu}}^m \right\rangle%
+ \left\langle \hat{Q}^\var(t),\hat{\vth}^m\right\rangle.
\end{split}
\end{equation}
Notice that we shall use the notation $|\cdot|_0$ for a (vector or scalar) $L^2$ norm. The same applies for $\|\cdot\|_1$ to denote a $H^1$ norm. Integrating in $[0,t]$, taking into account (\ref{initial_conditions}), the coercivity of $a^{V,\var}$, $a^{S,\var}$, integrating by parts the term in $\hat{\Gamma}_N^\var$ and using Korn's inequality we get
\begin{equation}\label{Galerkin-Eq5}
\begin{split}
&\hat{\rho}^\var \left| \dot{\hat{\bu}}^m(t) \right|^2_{0}%
+C \left\|{\hat{\bu}}^m(t)\right\|^2_{V}%
+\hat{\beta}^\var\left| \hat{\vth}^m (t)\right|^2_{0}%
+\hat{k}\tilde{C}\int_0^t\left\|{\hat{\vth}}^m(s)\right\|^2_{S}ds\\%
&\quad\leq  \int_0^t \left\{\left|\hat{\fb}^\var (s)\right|_{0}\left|\dot{\hat{\bu}}^m(s)\right|_{0}  + \left|\dot{\hat{\bh}}^\var (s)\right|_{0,\hat{\Gamma}_N^\var}\left|\hat{\bu}^m(s)\right|_{0,\hat{\Gamma}_N^\var}  +
\left|\hat{q}^\var(s)\right|_{0}\left|\hat{\vth}^m(s)\right|_{0}\right\}ds.
\end{split}
\end{equation}
Above and in what follows, $C, \tilde{C}$ denote positive constants whose specific value may change from line to line, only depending on data. Next, applying Young's inequality to each term in the right side in \eqref{Galerkin-Eq5} and the continuity of the trace operator, yields that
\begin{equation}\label{Galerkin-Eq6}
\begin{split}
&\left| \dot{\hat{\bu}}^m (t)\right|^2_{0}
+ \left\|{\hat{\bu}}^m(t)\right\|^2_{V}+
  \left| \hat{\vth}^m (t)\right|^2_{0}  + \int_0^t\left\|{\hat{\vth}}^m(s)\right\|^2_{S}ds\\
&\quad\leq C(\hat{\fb}^\var,\dot{\hat{\bh}}^\var,\hat{q}^\var) + \tilde{C} \int_0^t \left\{\left|\dot{\hat{\bu}}^m(s)\right|^2_{0}  + \left\|{\hat{\bu}}^m(s)\right\|^2_{1}+ \left|\hat{\vth}^m(s)\right|^2_{0}\right\}ds,
\end{split}
\end{equation}
%
which, applying Gronwall's Lemma, gives
\begin{equation}\label{Galerkin-cotas}
\left| \dot{\hat{\bu}}^m(t) \right|^2_{0}
+ \left\|{\hat{\bu}}^m(t)\right\|^2_{V}+   \left| \hat{\vth}^m (t)\right|^2_{0}
\leq C(\hat{\fb}^\var,\dot{\hat{\bh}}^\var,\hat{q}^\var) + e^{\tilde{C} T} , \, \, \forall \, m,
\end{equation}
from where,
\begin{equation*}
\dot{\hat{\bu}}^m \,\in \, L^\infty(0,T,[L^2(\hat{\Omega}^\var)]^3),\quad
\hat{\vth}^m \,\in \, L^\infty(0,T,L^2(\hat{\Omega}^\var)),\quad
\hat{\bu}^m \,\in \, L^\infty(0,T,V(\hat{\Omega}^\var)).
\end{equation*}
Further, going back to \eqref{Galerkin-Eq6}, we have
\begin{equation*}
 \hat{\vth}^m \,\in \, L^2(0,T,S(\hat{\Omega}^\var)).
\end{equation*}
%
Note that all the estimates are independent of $m$. Then

\begin{equation}\label{Exist-1}
 \left\{ \hat{\bu}^m\right\}_m \text{is a  bounded subset of }\, L^\infty(0,T,V(\hat{\Omega}^\var)) ,\end{equation}
\begin{equation}\label{Exist-2}
\left\{{\dot{\hat{\bu}}}^m \right\}_m  \text{is a  bounded subset of }\,\, \,  L^\infty(0,T,[L^2(\hat{\Omega}^\var)]^3),
\end{equation}
\begin{equation}\label{Exist-3}
 \left\{ \hat{\vth}^m\right\}_m \text{is a  bounded subset of }\,  \,\, L^\infty(0,T,L^2(\hat{\Omega}^\var))   \text{ and } \,   L^2(0,T;S(\hat{\Omega}^\var)).
\end{equation}
We now add equations (\ref{Galerkin-Eq1}) and (\ref{Galerkin-Eq2}) and write the result at times $t+h$, with $h>0$ and $0\le t\le T-h$, then subtract the resulting equations to get:
\begin{align*}
  &\hat{\rho}^\var \left\langle \ddot{\hat{u}}^m_i(t+h)-\ddot{\hat{u}}^m_i(t),\hat{v}_i^m \right\rangle%
+a^{V,\var}(\hat{\bu}^m(t+h)-\hat{\bu}^m(t),\hat{\bv}^m)%
-c^\var(\hat{\vartheta}^m(t+h)-\hat{\vartheta}^m(t),\hat{\bv}^m)\\
&\qquad
+\hat{\beta}^\var \left\langle\dot{\hat{\vartheta}}^m(t+h)-\dot{\hat{\vartheta}}^m(t), \hat{\varphi}^m\right\rangle\\%
&\qquad+a^{S,\var}(\hat{\vartheta}^m(t+h)-\hat{\vartheta}^m(t),\hat{\varphi}^m)%
+c^\var(\hat{\varphi}^m,\dot{\hat{\bu}}^m(t+h)-\dot{\hat{\bu}}^m(t))
\\%
&\qquad=\int_{\hat{\Omega}^\var}(\hat{f}^{i,\varepsilon}(t+h)-\hat{f}^{i,\varepsilon}(t))\hat{v}_i^m\,d\hat{x}^\varepsilon%
+\int_{\hat{\Gamma}_N^\var}(\hat{h}^{i,\varepsilon}(t+h)-\hat{h}^{i,\varepsilon}(t))\hat{v}_i^m\,d\hat{\Gamma}^\var +\int_{\hat{\Omega}^\var}(\hat{q}^\var(t+h)-\hat{q}^\var(t))\hat{\varphi}^m\,d\hat{x}^\varepsilon\\%
&\qquad \forall \hat{\bv}^m\in V_m,\ \forall \hat{\varphi}^m\in S_m,\aes.
\end{align*}
Next we take $\hat{\bv}^m=\dot{\hat{\bu}}^m(t+h)-\dot{\hat{\bu}}^m(t)\in V_m$ and $\hat{\varphi}^m=\hat{\vartheta}^m(t+h)-\hat{\vartheta}^m(t)\in S_m$ to obtain
\begin{align*}
  &\hat{\rho}^\var \left\langle \ddot{\hat{u}}^m_i(t+h)-\ddot{\hat{u}}^m_i(t),\dot{\hat{u}}^m_i(t+h)-\dot{\hat{u}}_i^m(t) \right\rangle%
+a^{V,\var}(\hat{\bu}^m(t+h)-\hat{\bu}^m(t),\dot{\hat{\bu}}^m(t+h)-\dot{\hat{\bu}}^m(t))\\%
&\qquad+\hat{\beta}^\var \left\langle\dot{\hat{\vartheta}}^m(t+h)-\dot{\hat{\vartheta}}^m(t), \hat{\vartheta}^m(t+h)-\hat{\vartheta}^m(t)\right\rangle %
+a^{S,\var}(\hat{\vartheta}^m(t+h)-\hat{\vartheta}^m(t),\hat{\vartheta}^m(t+h)-\hat{\vartheta}^m(t))\\%
&\qquad=\int_{\hat{\Omega}^\var}(\hat{f}^{i,\varepsilon}(t+h)-\hat{f}^{i,\varepsilon}(t))(\dot{\hat{u}}^m_i(t+h)-\dot{\hat{u}}_i^m(t))\,d\hat{x}^\varepsilon%
+\int_{\hat{\Gamma}_N^\var}(\hat{h}^{i,\varepsilon}(t+h)-\hat{h}^{i,\varepsilon}(t))(\dot{\hat{u}}^m_i(t+h)-\dot{\hat{u}}_i^m(t))\,d\hat{\Gamma}^\var \\&\qquad +\int_{\hat{\Omega}^\var}(\hat{q}^\varepsilon(t+h)-\hat{q}^\varepsilon(t))(\hat{\vartheta}^m(t+h)-\hat{\vartheta}^m(t))\,d\hat{x}^\varepsilon.
\end{align*}
Integrating in time in $[0,t]$ we get:
\begin{align*}
  &\frac{1}{2}\hat{\rho}^\var\left|\dot{\hat{\bu}}^m(t+h)-\dot{\hat{\bu}}^m(t)\right|_0^2%
-\frac{1}{2}\hat{\rho}^\var \left|\dot{\hat{\bu}}^m(h)-\dot{\hat{\bu}}^m(0)\right|_0^2 \\ &\qquad+\frac{1}{2}a^{V,\var}(\hat{\bu}^m(t+h)-\hat{\bu}^m(t),\hat{\bu}^m(t+h)-\hat{\bu}^m(t))%
-\frac{1}{2}a^{V,\var}(\hat{\bu}^m(h)-\hat{\bu}^m(0),\hat{\bu}^m(h)-\hat{\bu}^m(0))\\
&\qquad+\frac{1}{2}\int_{\hat{\Omega}^\var}\hat{\beta}^\var (\hat{\vartheta}^m(t+h)-\hat{\vartheta}^m(t))^2 d\hat{x}^\varepsilon-\frac{1}{2}\int_{\hat{\Omega}^\var}\hat{\beta}^\var (\hat{\vartheta}^m(h)-\hat{\vartheta}^m(0))^2 d\hat{x}^\varepsilon\\%
&\qquad+\int_{0}^{t}a^{S,\var}(\hat{\vartheta}^m(r+h)-\hat{\vartheta}^m(r),\hat{\vartheta}^m(r+h)-\hat{\vartheta}^m(r))dr\\%
&\qquad
=\int_{0}^{t}\int_{\hat{\Omega}^\var}(\hat{f}^{i,\varepsilon}(r+h)-\hat{f}^{i,\varepsilon}(r)) (\dot{\hat{u}}^m_i(r+h)-\dot{\hat{u}}_i^m(r))\,d\hat{x}^\varepsilon dr \\&\qquad%
+\int_{0}^{t}\int_{\hat{\Gamma}_N^\var}(\hat{h}^{i,\varepsilon}(r+h)-\hat{h}^{i,\varepsilon}(r)) (\dot{\hat{u}}^m_i(r+h)-\dot{\hat{u}}_i^m(r))\,d\hat{\Gamma}^\var dr \\&\qquad +\int_{0}^{t}\int_{\hat{\Omega}^\var}(\hat{q}^\varepsilon(r+h)-\hat{q}^\varepsilon(r)) (\hat{\vartheta}^m(r+h)-\hat{\vartheta}^m(r))\,d\hat{x}^\varepsilon dr.
\end{align*}
Now, dividing the equation by $h^2$ and having in mind \eqref{Galerkin-cotas},
we can take limits when $h\rightarrow 0^+$ to have
\begin{align}
  &\frac{1}{2}\hat{\rho}^\var\left|\ddot{\hat{\bu}}^m(t)\right|^2_{0} -\frac{1}{2}\hat{\rho}^\var\left|\ddot{\hat{\bu}}^m(0) \right|^2_{0}
  +\frac{1}{2}a^{V,\var}(\dot{\hat{\bu}}^m(t),\dot{\hat{\bu}}^m(t))
  -\frac{1}{2}a^{V,\var}(\dot{\hat{\bu}}^m(0),\dot{\hat{\bu}}^m(0))\nonumber\\
&\qquad+\frac{1}{2}\int_{\hat{\Omega}^\var}\hat{\beta}^\var (\dot{\hat{\vartheta}}^m(t))^2 d\hat{x}^\varepsilon -\frac{1}{2}\int_{\hat{\Omega}^\var}\hat{\beta}^\var (\dot{\hat{\vartheta}}^m(0))^2 d\hat{x}^\varepsilon%
+\int_{0}^{t}a^{S,\var}(\dot{\hat{\vartheta}}^m(r),\dot{\hat{\vartheta}}^m(r))dr\nonumber\\%
&\qquad
=\int_{0}^{t}\int_{\hat{\Omega}^\var}\dot{\hat{f}}^{i,\varepsilon}(r)\ddot{\hat{u}}^m_i(r)\, d\hat{x}^\varepsilon dr%
+\int_{0}^{t}\int_{\hat{\Gamma}_N^\var}\dot{\hat{h}}^{i,\var}(r)\ddot{\hat{u}}^m_i(r) \,d\hat{\Gamma}^\var dr  +\int_{0}^{t}\int_{\hat{\Omega}^\var}\dot{\hat{q}}^\var(r)\dot{\hat{\vartheta}}^m(r)\,d\hat{x}^\var dr.\label{paraluego}
\end{align}
Integrating by parts the term on $\hat{\Gamma}_N^\var$ above and applying Young's inequality, we get
\begin{align}
  &\hat{\rho}^\var |\ddot{\hat{\bu}}^m(t)|^2_{0}-\hat{\rho}^\var |\ddot{\hat{\bu}}^m(0)|^2_{0}
+||\dot{\hat{\bu}}^m(t)||^2_{V}
+\hat{\beta}^\var |\dot{\hat{\vartheta}}^m(t)|_0^2 -\hat{\beta}^\var |\dot{\hat{\vartheta}}^m(0)|_0^2
+\int_{0}^{t}\|\dot{\hat{\vartheta}}^m(r)\|_S^2 dr\nonumber\\%
&\qquad\le \tilde{C}(\dot{\hat{\fb}}^\var,\ddot{\hat{\bh}}^\var,\dot{\hat{q}}^\var)+ C\int_{0}^{t}\left\{|\ddot{\hat{\bu}}^m(r)|^2_{0}+||\dot{\hat{\bu}}^m(r)||^2_{1}+ |\dot{\hat{\vartheta}}^m(r)|_0^2\right\}dr.\label{cotas_FG}
\end{align}
In order to obtain bounds for $|\ddot{\hat{\bu}}^m(0)|^2_{0}$ and $|\dot{\hat{\vartheta}}^m(0)|_0^2$ we first notice that equations (\ref{Galerkin-Eq1}) and (\ref{Galerkin-Eq2}) hold for $t=0$ due to the compatibility required between initial and boundary conditions. Therefore, taking $t=0$ and $\hat{\bv}^m=\ddot{\hat{\bu}}^m(0)\in V_m$ in (\ref{Galerkin-Eq1}) and $\hat{\varphi}^m=\dot{\hat{\vartheta}}^m(0)\in S_m$ in (\ref{Galerkin-Eq2}), taking into account the initial conditions, and using Young's inequality, we obtain 
\begin{align*}
&\hat{\rho}^\var|\ddot{\hat{\bu}}^m(0)|_0^2=\int_{\hat{\Omega}^\var}\hat{f}^{i,\varepsilon}(0)\ddot{\hat{u}}^m_i(0)\,d\hat{x}^\varepsilon%
+\int_{\hat{\Gamma}_N^\var}\hat{h}^{i,\varepsilon}(0)\ddot{\hat{u}}^m_i(0)\,d\hat{\Gamma}^\var\le \frac{1}{\delta}C+\delta|\ddot{\hat{\bu}}^m(0)|_0^2,\\
&\hat{\beta}^\var| \dot{\hat{\vartheta}}^m(0)|_0^2 =
\int_{\hat{\Omega}^\var}\hat{q}^\varepsilon(0)\dot{\hat{\vartheta}}^m(0)\,d\hat{x}^\varepsilon   \ \le \frac{1}{\tilde{\delta}}\tilde{C}+\tilde{\delta}| \dot{\hat{\vartheta}}^m(0)|_0^2,
\end{align*}
where $\delta$, and $\tilde{\delta}$ are sufficiently small positive constants. Next, applying Korn's inequality and Gronwall's lemma in (\ref{cotas_FG}) we find
\begin{equation*}
  |\ddot{\hat{\bu}}^m(t)|^2_{0}+||\dot{\hat{\bu}}^m(t)||^2_{V}+|\dot{\hat{\vartheta}}^m(t)|_0^2\le C.
\end{equation*}
Again, all the estimates are independent of $m$. Then,

\begin{align}\label{Exist-8}
 &\left\{ \dot{\hat{\bu}}^m\right\}_m \text{is a  bounded subset of }\, L^\infty(0,T,V(\hat{\Omega}^\var)) ,\\
\label{Exist-9}
&\left\{{\ddot{\hat{\bu}}}^m \right\}_m  \text{is a  bounded subset of }\,\, \,  L^\infty(0,T,[L^2(\hat{\Omega}^\var)]^3),
\\\label{Exist-10}
 &\left\{ \dot{\hat{\vth}}^m\right\}_m \text{is a  bounded subset of }\,  \,\, L^\infty(0,T,L^2(\hat{\Omega}^\var)).
\end{align}
Observe that \eqref{Exist-1}--\eqref{Exist-3} and \eqref{Exist-8}--\eqref{Exist-10} imply that there exists  subsequences of $\hat{\bu}^m$ and $\hat{\vartheta}^m$, also denoted by $\hat{\bu}^m$ and $\hat{\vartheta}^m$, and there exist elements $\hat{\bu}^\var$, $\dot{\hat{\bu}}^\var$, $\ddot{\hat{\bu}}^\var$, $\hat{\vartheta}^\var$, $\dot{\hat{\vartheta}}^\var$ 
such that
\begin{align}
&\hat{\bu}^m \stackrel[m\to\infty]{ }{\xrightharpoonup{\ \ *\ \ }}\hat{\bu}^\var \quad {\rm in }\, L^\infty(0,T;V(\hat{\Omega}^\var)), \label{Exist-5}\\
&\dot{\hat{\bu}}^m\stackrel[m\to\infty]{ }{\xrightharpoonup{\ \ *\ \ }}\dot{\hat{\bu}}^\var \quad {\rm in }\, L^\infty(0,T;[L^2(\hat{\Omega}^\var)]^3)\cap L^\infty(0,T;V(\hat{\Omega}^\var)),\label{Exist-6}\\
&\ddot{\hat{\bu}}^m\stackrel[m\to\infty]{ }{\xrightharpoonup{\ \ *\ \ }}\ddot{\hat{\bu}}^\var \quad {\rm in }\, L^\infty(0,T;[L^2(\hat{\Omega}^\var)]^3),\label{Exist-6_ac}\\
&\hat{\vartheta}^m \stackrel[m\to\infty]{ }{\xrightharpoonup{\ \ *\ \ }} \hat{\vartheta}^\var \quad {\rm in }\, L^\infty(0,T;L^2(\hat{\Omega}^\var))\cap L^\infty(0,T;S(\hat{\Omega}^\var)), \label{Exist-5-bis}\\
&\dot{\hat{\vartheta}}^m \stackrel[m\to\infty]{ }{\xrightharpoonup{\ \ *\ \ }} \dot{\hat{\vartheta}}^\var \quad {\rm in }\, L^\infty(0,T;L^2(\hat{\Omega}^\var)). \label{Exist-5-velt}\\
\end{align}
Consider now $\hat{\bv}^m=\hat{\bw}_j$ and $\hat{\varphi}^m=\hat{s}_i$ in equations \eqref{Galerkin-Eq1} and \eqref{Galerkin-Eq2} fixed:

\begin{equation}\label{Galerkin-Eq8}
\hat{\rho}^\var \left\langle \ddot{\hat{\bu}}^m,\hat{\bw}_j \right\rangle
+ a^{V,\var}({\hat{\bu}}^m,\hat{\bw}_j)-c^\var(\hat{\vartheta}^m,\hat{\bw}_j)
=\left\langle \hat{J}^\var(t),\hat{\bw}_j \right\rangle,
	\end{equation}
	
\begin{equation}\label{Galerkin-Eq9}
\left\langle\dot{\hat{\vth}}^m,\hat{s}_i \right\rangle
+ a^{S,\var}({\hat{\vth}}^m,\hat{s}_i)
+ c^\var(\hat{s}_i,\dot{\hat{\bu}}^m)
=\left\langle \hat{Q}^\var(t),\hat{s}_i\right\rangle.
\end{equation}
 Observe that \eqref{Exist-5} and \eqref{Exist-6} imply that
\begin{equation*}
a^{V,\var}({\hat{\bu}}^m,\hat{\bw}_j)\rightarrow a^{V,\var}({\hat{\bu}}^\var,\hat{\bw}_j) \, \, \text{ and }
 c^\var(\hat{s}_i,\dot{\hat{\bu}}^m) \rightarrow c^\var(\hat{s}_i,\dot{\hat{\bu}}^\var)  \, \,\text{ in }\, L^\infty(0,T).
\end{equation*}
%
%
Analogously, from \eqref{Exist-5-bis} we can state that
\begin{equation*}
 a^{S,\var}(\hat{\vth}^m ,s_i) \rightarrow a^{S,\var}(\hat{\vth}^\var,s_i)  \, \, \text{ and }
  c^\var(\hat{\vth}^m ,\hat{\bw}_j) \rightarrow c^\var(\hat{\vth}^\var,\hat{\bw}_j)  \, \, \text{ in }\, L^\infty(0,T).
 \end{equation*}
Now, from (\ref{Exist-6_ac}) and (\ref{Exist-5-velt}) we have:
\begin{equation*}
\left\langle\ddot{\hat{\bu}}^m,\hat{\bw}_j\right\rangle=(\ddot{\hat{\bu}}^m,\hat{\bw}_j) \rightarrow (\ddot{\hat{\bu}}^\var,\hat{\bw}_j)  \, \, \text{ and }  \left\langle\dot{\vartheta^m},s_i\right\rangle=(\dot{\vartheta^m},s_i)\rightarrow (\dot{\vartheta}^\var,s_i)\qquad \, \text{ in }  L^\infty(0,T).
\end{equation*}
%
%
%
%
%
%
Then, we can take $m\rightarrow \infty$ in \eqref{Galerkin-Eq8}--\eqref{Galerkin-Eq9} obtaining that
\begin{align}
\hat{\rho}^\var (\ddot{\hat{\bu}}^\var,\hat{\bw}_j)
+ a^{V,\var}({\hat{\bu}^\var},\hat{\bw}_j)-c^\var(\hat{\vartheta}^\var,\hat{\bw}_j)
&=\left\langle \hat{J}^\var(t),\hat{\bw}_j \right\rangle, &\text{ in } \ 
L^\infty(0,T),\, \forall \, j\geq 1,\label{Galerkin-Eq8-lim}\\
(\dot{\hat{\vth}}^\var,\hat{s}_i)
+ a^{S,\var}(\hat{\vth}^\var,\hat{s}_i)
+ c^\var(\dot{\hat{\bu}}^\var,\hat{s}_i)
&=\left\langle \hat{Q}^\var(t),\hat{s}_i\right\rangle, &\text{ in } \, L^\infty(0,T),\, \forall \, i\geq 1.\label{Galerkin-Eq9-lim}
\end{align}
%

Besides, since the initial conditions (\ref{initial_conditions}) are null, it is trivial that, when $m\to\infty$, the limit functions have null initial conditions as well,  which completes the proof for the existence and regularity of the solutions. We focus now on proving the uniqueness.%

%
%
%
%
%

Let us assume that there exist two solutions $\{\hat{\bu}^{\var,1},\hat{\vartheta}^{\var,1}\}$ and $\{\hat{\bu}^{\var,2},\hat{\vartheta}^{\var,2}\}$ for Problem \ref{problema_cart}. Let us define $\bw^\var=\hat{\bu}^{\var,1}-\hat{\bu}^{\var,2}$ and $\phi^\var=\hat{\vartheta}^{\var,1}-\hat{\vartheta}^{\var,2}.$
Now, we consider equations (\ref{P2.2a})--(\ref{P2.2b}) at time $t$ for $\{\hat{\bu}^{\var,i},\hat{\vartheta}^{\var,i}\}$, take as test function $\hat{\bv}^\var=\dot{\bw}^\var(t)$ and $\hat{\varphi}^\var=\phi^\var(t)$ for both $i=1$ and $i=2$, and subtract the resulting equations to find:
\begin{align*}
&\int_{\hat{\Omega}^\var}\hat{\rho}^\var \ddot{\bw}^\var(t) \dot{\bw}^\var(t) d\hat{x}^\varepsilon%
+a^{V,\var}(\bw^\var(t),\dot{\bw}^\var(t))%
-c^\var(\phi^\var(t),\dot{\bw}^\var(t))=0\\
&\int_{\hat{\Omega}^\var}\hat{\beta}^\var \dot{\phi}^\var(t) \phi^\var(t) d\hat{x}^\varepsilon%
+a^{S,\var}(\phi^\var(t),\phi^\var(t))%
+c^\var(\phi^\var(t),\dot{\bw}^\var(t))=0,\quad\aes.
\end{align*}
Adding these last two equations we have
\begin{align*}
&\int_{\hat{\Omega}^\var}\hat{\rho}^\var \ddot{\bw}^\var(t) \dot{\bw}^\var(t) d\hat{x}^\varepsilon%
+a^{V,\var}(\bw^\var(t),\dot{\bw}^\var(t))%
\nonumber\\%
&\qquad + \int_{\hat{\Omega}^\var}\hat{\beta}^\var \dot{\phi}^\var(t) \phi^\var(t) d\hat{x}^\varepsilon%
+a^{S,\var}(\phi^\var(t),\phi^\var(t))=0,\quad\aes.
\end{align*}
Integrating in $[0,t]$, and taking into account the initial conditions we obtain:
\begin{align*}
&\hat{\rho}^\var |\dot{\bw}^\var(t)|_0^2 +||\bw^\var(t)||_V^2+\hat{\beta}^\var |\phi^\var(t)|_0^2
+\int_0^t a^{S,\var}(\phi^\var(r),\phi^\var(r))dr
=0,\quad\aes,
\end{align*}
from where one easily deduce that $\bw^\var=\bzero$ and $\phi^\var
=0.$
\qed

\section{A three-dimensional dynamic problem for thermoelastic shells}\label{seccion_dominio_ind}

In this section we consider the particular case when the deformable body is, in fact, a shell. We first provide key notations and some preliminary results in a summarised form. The interested reader can consult \cite{Ciarlet4b} and \cite{ContactShell} for a more detailed exposition.

Let $\omega$ be a bounded domain of $\mathbb{R}^2$, with a Lipschitz-continuous boundary $\gamma=\partial\omega$. Let $\by=(y_\alpha)$ be a generic point of  its closure $\bar{\omega}$ and let $\d_\alpha$ denote the partial derivative with respect to $y_\alpha$. Above and in what follows, Greek indices take their values in the set  $\{1,2\}$, whereas Latin indices do it in the set $\{1,2,3\}$. We will use summation convention on repeated indices. Let $\btheta\in\mathcal{C}^2(\bar{\omega};\mathbb{R}^3)$ be an injective mapping such that the two vectors $\ba_\alpha(\by):= \d_\alpha \btheta(\by)$ are linearly independent. These vectors form the covariant basis of the tangent plane to the surface $S:=\btheta(\bar{\omega})$ at the point $\btheta(\by).$ We also consider the two vectors $\ba^\alpha(\by)$ of the same tangent plane defined by the relations $\ba^\alpha(\by)\cdot \ba_\beta(\by)=\delta_\beta^\alpha$, that constitute its contravariant basis. We define
$$
\ba_3(\by)=\ba^3(\by):=\frac{\ba_1(\by)\wedge \ba_2(\by)}{| \ba_1(\by)\wedge \ba_2(\by)|},
$$
the unit normal vector to $S$ at the point $\btheta(\by)$, where $\wedge$ denotes vector product in $\mathbb{R}^3$. We can define the first fundamental form, given as metric tensor, in covariant or contravariant components, respectively, by $a_{\alpha\beta}:=\ba_\alpha\cdot \ba_\beta$, $a^{\alpha\beta}:=\ba^\alpha\cdot \ba^\beta$. The second fundamental form, given as curvature tensor, in covariant or mixed components, respectively, is given by $b_{\alpha\beta}:=\ba^3 \cdot \d_\beta \ba_\alpha$, $b_{\alpha}^\beta:=a^{\beta\sigma}\cdot b_{\sigma\alpha}$, and the Christoffel symbols of the surface $S$ as
$\Gamma^\sigma_{\alpha\beta}:=\ba^\sigma\cdot \d_\beta \ba_\alpha$. The area element along $S$ is $\sqrt{a}dy$ where $a:=\det (a_{\alpha\beta})$.

We define the three-dimensional domain $\Omega^\varepsilon:=\omega \times (-\varepsilon, \varepsilon)$ and its boundary $\Gamae=\partial\Omega^\var$. We also define  the following parts of the boundary, $\Gamma^\varepsilon_N:=\omega\times \{\varepsilon\}$ (it could also be the lower face or the union of both), 
$\Gamma_0^\varepsilon:=\gamma_0\times[-\varepsilon,\varepsilon]$, where $\gamma_0\subseteq\gamma$. Let $\bx^\varepsilon=(x_i^\varepsilon)$ be a generic point of $\bar{\Omega}^\varepsilon$ and let $\d_i^\var$ denote the partial derivative with respect to $x_i^\varepsilon$. Note that $x_\alpha^\varepsilon=y_\alpha$ and $\d_\alpha^\varepsilon =\d_\alpha$. Let $\bTheta:\bar{\Omega}^\varepsilon\rightarrow \mathbb{R}^3$ be the mapping defined by
\begin{equation} \label{bTheta}
\bTheta(\bx^\varepsilon):=\btheta(\by) + x_3^\varepsilon \ba_3(\by) \ \forall \bx^\varepsilon=(\by,x_3^\varepsilon)=(y_1,y_2,x_3^\varepsilon)\in\bar{\Omega}^\varepsilon.
\end{equation}
By identifying $\hat{\Omega}^\var=\bTheta(\Omega^\var)$, $\hat{\Gamma}^\var=\bTheta(\Gamma^\var)$, $\hat{\Gamma}^\var_0=\bTheta(\Gamma^\var_0)$, etc. we cast this setting into the more general three dimensional framework of the preceding section, as a particular case. Further,
in \cite[Th.  3.1-1]{Ciarlet4b} it is shown that if the injective mapping $\btheta:\bar{\omega}\rightarrow\mathbb{R}^3$ is smooth enough, the mapping $\bTheta:\bar{\Omega}^\var\rightarrow\mathbb{R}^3$ is also injective for $\var>0$ small enough and the vectors $\bg_i^\varepsilon(\bx^\varepsilon):=\d_i^\varepsilon\bTheta(\bx^\varepsilon)$ are linearly independent.
Therefore, the three vectors $\bg_i^\varepsilon(\bx^\varepsilon)$ form the covariant basis at the point $\bTheta(\bx^\varepsilon)$ and $\bg^{i,\varepsilon}(\bx^\varepsilon) $ defined by the relations
$\bg^{i,\varepsilon}\cdot \bg_j^\varepsilon=\delta_j^i$ form the contravariant basis at the point $\bTheta(\bx^\varepsilon)$. The covariant and contravariant components of the metric tensor are defined, respectively, as $g_{ij}^\varepsilon:=\bg_i^\varepsilon \cdot \bg_j^\varepsilon$, $g^{ij,\varepsilon}:=\bg^{i,\varepsilon} \cdot \bg^{j,\varepsilon}$, and Christoffel symbols as
$\Gamma^{p,\varepsilon}_{ij}:=\bg^{p,\varepsilon}\cdot\d_i^\varepsilon \bg_j^\varepsilon$.
The volume element in the set $\bTheta(\bar{\Omega}^\varepsilon)$ is $\sqrt{g^\varepsilon}dx^\var$ and the surface element in $\bTheta(\Gamma^\varepsilon)$ is $\sqrt{g^\varepsilon}d\Gamae$  where
$g^\varepsilon:=\det (g^\varepsilon_{ij})$.

We now define the corresponding contravariant components in curvilinear coordinates for the applied forces densities:
$$
\hat{f}^{i,\var}(\hat{\bx}^\var)\hat{\be}_i\,d\hat{x}^\var=:f^{i,\var}(\bx^\var)\bg_i^\var(\bx^\var)\sqrt{g^\var(\bx^\var)}\,dx^\var,\quad%
\hat{h}^{i,\var}(\hat{\bx}^\var)\hat{\be}_id\hat{\Gamma}^\var=:h^{i,\var}(\bx^\var)\bg_i^\var(\bx^\var)\sqrt{g^\var(\bx^\var)}d\Gamma^\var, $$
and the covariant components in curvilinear coordinates for the displacements field:
$$
\hat{\bu}^\var(\hat{\bx}^\var)=\hat{u}^\var_i(\hat{\bx}^\var)\hat{\be}^i=:u_i^\var(\bx^\var)\bg^{i,\var}(\bx^\var),\ {\rm with}\ \hat{\bx}^\var=\bTheta(\bx^\var).
$$
\begin{remark}
Notice that forces above depend also on the time variable $t\in[0,T]$, but we decided to keep it implicit for the sake of readiness, since the subject of the change of variable is the spatial component. The same comment applies in a number of situations below.
\end{remark}
We also define $\vartheta^\var(\bx^\var):=\hat{\vartheta}^\var(\hat{\bx}^\var)$ and $q^\var(\bx^\var):=\hat{q}^\var(\hat{\bx}^\var)$.

Let us define the spaces,
\begin{align*} 
&V(\Omega^\varepsilon)=\{\bv^\varepsilon=(v_i^\varepsilon)\in [H^1(\Omega^\varepsilon)]^3; \bv^\varepsilon=\mathbf{\bcero} \ {\rm on} \ \Gamma_0^\varepsilon\},
&{S}(\Omega^\varepsilon)=\{\varphi^\varepsilon\in H^1(\Omega^\varepsilon); \varphi^\varepsilon=0 \ {\rm on} \ \Gamma_0^\varepsilon\}.
\end{align*}
Both are real Hilbert spaces with the induced inner product of $[H^1(\Omega^\var)]^d$, $d\in\{1,3\}$. The corresponding norm is denoted by $\left\|\cdot\right\|_{1,\Omega^\var}$ in both cases, since no confusion is possible.  With these definitions it is straightforward to derive from the Problem \ref{problema_cart} the following variational problem (see \cite{Ciarlet4b} for the case in linear elasticity and use similar arguments):
\begin{problem}\label{problema_eps}
Find a pair $t\mapsto(\bu^\varepsilon(\bx^\var,t),\vartheta^\varepsilon(\bx^\var,t))$ of $[0,T]\to V(\Omega^\var)\times {S}(\Omega^\var)$ verifying
\begin{align*}
&\int_{\Omega^\var}\rho^\var (\ddot{u}_\alpha^\varepsilon g^{\alpha\beta,\var}v_\beta^\varepsilon+\ddot{u}_3^\varepsilon v_3^\varepsilon)\sqrt{g^\varepsilon}\,d{x}^\varepsilon%
+\int_{\Omega^\varepsilon}A^{ijkl,\varepsilon}e^\varepsilon_{k||l}(\bu^\varepsilon)e^\varepsilon_{i||j}(\bv^\varepsilon)\sqrt{g^\varepsilon} dx^\varepsilon\\%
&\qquad-\int_{\Omega^\var}\alpha_T^\var(3\lambda^\var+2\mu^\var)\vartheta^\var(\eab^\var(\bv^\varepsilon)g^{\alpha\beta,\var}+\edtres^\var(\bv^\var))\sqrt{g^\varepsilon} dx^\varepsilon%
\\%
&\qquad= \int_{\Omega^\varepsilon} f^{i,\varepsilon}v_i^\varepsilon\sqrt{g^\varepsilon} dx^\varepsilon%
  +\int_{\Gamma_N^\varepsilon} h^{i,\varepsilon}v_i^\varepsilon\sqrt{g^\varepsilon} d\Gamma^\varepsilon  \quad \forall \bv^\varepsilon\in V(\Omega^\varepsilon),\aes,\\
&\int_{\Omega^\var}\beta^\var\dot{\vartheta}^\varepsilon \varphi^\varepsilon \sqrt{g^\varepsilon} dx^\varepsilon%
+\int_{\Omega^\var}k^\varepsilon(\partial_\alpha^\varepsilon \vartheta^\varepsilon g^{\alpha\beta,\var}\partial_\beta^\varepsilon\varphi^\varepsilon+\partial_3^\varepsilon \vartheta^\varepsilon\partial_3^\varepsilon\varphi^\varepsilon)\sqrt{g^\varepsilon} dx^\varepsilon\\%
&\qquad+\int_{\Omega^\var}\alpha_T^\var(3\lambda^\var+2\mu^\var)\varphi^\var(\eab^\var(\dot{\bu}^\varepsilon)g^{\alpha\beta,\var}+\edtres^\var(\dot{\bu}^\var))\sqrt{g^\varepsilon} dx^\varepsilon\\%
&\qquad=\int_{\Omega^\var}q^\varepsilon\varphi^\varepsilon\,\sqrt{g^\varepsilon} dx^\varepsilon   \quad \forall \varphi^\varepsilon\in {S}(\Omega^\var),\aes,
  \end{align*}
with $\dot{\bu}^\varepsilon(\cdot,0)=\bu^\varepsilon(\cdot,0)=\bzero$ and $\vartheta^\var(\cdot,0)=0$.
\end{problem}
Above, $A^{ijkl,\varepsilon}=A^{jikl,\varepsilon}=A^{klij,\varepsilon}\in{\mathcal{C}}^1(\bar{\Omega}^\var)$, defined by
\begin{equation}\label{TensorAeps}
A^{ijkl,\varepsilon}:= \lambda g^{ij,\varepsilon}g^{kl,\varepsilon} + \mu(g^{ik,\varepsilon}g^{jl,\varepsilon} + g^{il,\varepsilon}g^{jk,\varepsilon} ), 
\end{equation}
represent the contravariant components of the three-dimensional elasticity tensor, and the functions $e^\varepsilon_{i||j}(\bv^\var)=e^\varepsilon_{j||i}(\bv^\var)\in L^2(\Omega^\var)$ that represent the covariant components of the linearized change of metric tensor, or strain tensor, are defined  by
\begin{align*}
e^\varepsilon_{i||j}(\bv^\var):= 
\frac1{2}(\d^\varepsilon_jv^\varepsilon_i + \d^\varepsilon_iv^\varepsilon_j) - \Gamma^{p,\varepsilon}_{ij}v^\varepsilon_p,
\end{align*}
for all $\bv^\var\in [H^1(\Omega^\var)]^3$, where $\partial_i^\var$ denotes partial derivative with respect to $x_i^\var$. Note that the following simplifications are verified,
\begin{equation}
\Gamma^{3,\varepsilon}_{\alpha 3}=\Gamma^{p,\varepsilon}_{33}=0  \ \textrm{in} \ \bar{\Omega}^\varepsilon,%
\quad  A^{\alpha\beta\sigma 3,\varepsilon}=A^{\alpha 333,\varepsilon}=0 \ \textrm{in} \ \bar{\Omega}^\varepsilon,\label{tensor_terminos_nulos}
\end{equation}
as a consequence of the definition of $\bTheta$ in (\ref{bTheta}). The definitions of the fourth order tensor (\ref{TensorAeps}) imply that (see \cite[Theorem 1.8-1]{Ciarlet4b}) for $\var>0$ small enough, there exists a constant $C_e>0$, independent of $\var$, such that,
\begin{align} \label{elipticidadA}
 \sum_{i,j}|t_{ij}|^2\leq C_e A^{ijkl,\var}(\bx^\var)t_{kl}t_{ij},
\end{align}
for all $\bx^\var\in\bar{\Omega}^\var$ and all $\bt=(t_{ij})\in\mathbb{S}^3$ (vector space of $3\times 3$ real symmetric matrices).


\begin{remark}
We recall that the vector field $\bu^\varepsilon=(u_i^\varepsilon):{\Omega}^\varepsilon\times[0,T] \rightarrow \mathbb{R}^3$ solution of Problem \ref{problema_eps} has to be interpreted conveniently. The functions  $u_i^\varepsilon:\bar{\Omega}^\varepsilon\times[0,T] \rightarrow \mathbb{R}^3$ are the covariant, time dependent, components of the ``true" displacements field $\bUcal^\var:=u_i^\varepsilon \bg^{i,\varepsilon}:\bar{\Omega}^\varepsilon\times[0,T] \rightarrow \mathbb{R}^3$.
\end{remark}

For convenience, we consider a reference domain independent of the small parameter $\var$. Hence, let us define the three-dimensional domain $\Omega:=\omega \times (-1, 1) $ and  its boundary $\Gamma=\partial\Omega$. We also define the following parts of the boundary,
 \begin{align*}
 \Gamma_N:=\omega\times \{1\}, 
 \quad \Gamma_0:=\gamma_0 \times[-1,1].
 \end{align*}
 Let $\bx=(x_1,x_2,x_3)$ be a generic point in $\bar{\Omega}$ and we consider the notation $\d_i$ for the partial derivative with respect to $x_i$. We define the  projection map $\pi^\varepsilon: \bar{\Omega} \to\bar{\Omega}^\varepsilon,$ such that
 \begin{align*}
 \pi^\varepsilon(\bx)=\bx^\varepsilon=(x_i^\varepsilon)=(x_1^\var,x_2^\var,x_3^\var)=(x_1,x_2,\varepsilon x_3)\in \bar{\Omega}^\varepsilon,
 \end{align*}
 hence, $\d_\alpha^\varepsilon=\d_\alpha $  and $\d_3^\varepsilon=\frac1{\varepsilon}\d_3$. We consider the displacements related scaled unknown $\bu(\varepsilon)=(u_i(\varepsilon)):\bar{\Omega}\times[0,T]\to \mathbb{R}^3$ and the scaled vector fields $\bv=(v_i):\bar{\Omega}\to \mathbb{R}^3 $ defined as
 \begin{align*}
 u_i^\varepsilon(\bx^\varepsilon)=:u_i(\varepsilon)(\bx) \ \textrm{and} \ v_i^\varepsilon(\bx^\varepsilon)=:v_i(\bx) \ \forall \bx\in\bar{\Omega},\ \bx^\varepsilon=\pi^\varepsilon(\bx)\in \bar{\Omega}^\varepsilon.
 \end{align*}
Besides, we define the scaled temperature $\vartheta(\var):\bar{\Omega}\times[0,T]\to \mathbb{R}$ defined as
$$
\vartheta(\var)(\bx):=\vartheta^\var(\bx^\var)\quad \forall \bx\in\Omega, \ \textrm{where} \ \bx^\varepsilon=\pi^\varepsilon(\bx)\in \Omega^\varepsilon.
$$
For the sake of simplicity, from now on, we are going to assume that the different parameters of the problem (thermal conductivity, thermal dilatation, specific heat coefficient, mass density, Lamé coefficients) are all independent of $\var$.
Also, let the functions, $\Gamma_{ij}^{p,\varepsilon}, g^\varepsilon, A^{ijkl,\varepsilon}$ 
be associated with the functions $\Gamma_{ij}^p(\varepsilon),$ $ g(\varepsilon),$ $ A^{ijkl}(\varepsilon),$ defined by
$$
\Gamma_{ij}^p(\varepsilon)(\bx):=\Gamma_{ij}^{p,\varepsilon}(\bx^\varepsilon),\
g(\varepsilon)(\bx):=g^\varepsilon(\bx^\varepsilon),\
A^{ijkl}(\varepsilon)(\bx):=A^{ijkl,\varepsilon}(\bx^\varepsilon),
$$
for all $\bx\in\bar{\Omega}$, $\bx^\varepsilon=\pi^\varepsilon(\bx)\in\bar{\Omega}^\varepsilon$. For all $\bv=(v_i)\in [H^1(\Omega)]^3$, let there be associated the scaled linearized strains $(\eij(\var)(\bv))\in [L^2(\Omega)]^{3\times 3}_{sym}$, which we also denote as $(\eij(\var;\bv))$, defined by
\begin{align}
&\eab(\varepsilon;\bv):=\frac{1}{2}(\d_\beta v_\alpha + \d_\alpha v_\beta) - \Gamma_{\alpha\beta}^p(\varepsilon)v_p, \label{eab}\\
& \eatres(\varepsilon;\bv):=\frac{1}{2}(\frac{1}{\var}\d_3 v_\alpha + \d_\alpha v_3) - \Gamma_{\alpha 3}^p(\varepsilon)v_p,\label{eatres}\\
& \edtres(\varepsilon;\bv):=\frac1{\varepsilon}\d_3v_3. \label{edtres}
\end{align}
Note that with these definitions it is verified that
\begin{align*}
\eij^\var(\bv^\var)(\pi^\var(\bx))=\eij(\var;\bv)(\bx) \ \forall\bx\in\Omega.
\end{align*}
 \begin{remark} The functions $\Gamma_{ij}^p(\varepsilon), g(\varepsilon), A^{ijkl}(\varepsilon)$ converge in $\mathcal{C}^0(\bar{\Omega})$ when $\varepsilon$ tends to zero.
 \end{remark}
 \begin{remark}When we consider
 $\varepsilon=0$ the functions will be defined with respect to $\by\in\bar{\omega}$. Notice the singularities in (\ref{eatres}) and (\ref{edtres}) for that case. We shall distinguish the three-dimensional Christoffel symbols from the two-dimensional ones associated to $S$ by using  $\Gamma_{\alpha \beta}^\sigma(\varepsilon)$ and $ \Gamma_{\alpha\beta}^\sigma$, respectively.
 \end{remark}
In \cite[Theorem 3.3-2]{Ciarlet4b} we find an important result which shows that under suitable regularity conditions, take for example $\btheta\in\mathcal{C}^2(\bar{\omega};\mathbb{R}^3)$, there exists an $\var_0>0$ such that $A^{ijkl}(\var)$ is positive-definite, uniformly with respect to $\bx\in\bar{\Omega}$ and $\var$, provided that $0<\var\leq\var_0$. Further, the asymptotic behavior of $A^{ijkl}(\var)$ is detailed.
%
Indeed, it is satisfied that 
\begin{align*}
A^{ijkl}(\var)= A^{ijkl}(0) + {\mathcal{O} }(\var) \ \textrm{and} \ A^{\alpha\beta\sigma 3}(\var)=A^{\alpha 3 3 3}(\var)=0,
\end{align*}
for all $\var$, $0<\var \leq \var_0$, and
\begin{align}
&A^{\alpha\beta\sigma\tau}(0)=\lambda a^{\alpha\beta}a^{\sigma\tau} + \mu(a^{\alpha\sigma}a^{\beta\tau} + a^{\alpha\tau}a^{\beta\sigma}),\quad%
A^{\alpha\beta 3 3}(0)= \lambda a^{\alpha\beta},\label{eq51b1}\\
&A^{\alpha 3\sigma 3}(0)=\mu a^{\alpha\sigma},\quad%
A^{33 3 3}(0)= \lambda + 2\mu,\quad%
A^{\alpha\beta\sigma 3}(0)=A^{\alpha 333}(0)=0.\label{eq51b2}
\end{align}
Moreover, and related with \eqref{elipticidadA}, there exists a constant $C_e>0$, independent of the variables and $\var$, such that
  \begin{align} \label{elipticidadA_eps}
  \sum_{i,j}|t_{ij}|^2\leq C_e A^{ijkl}(\varepsilon)(\bx)t_{kl}t_{ij},
  \end{align}
 for all $\var$, $0<\var\leq\var_0$, for all $\bx\in\bar{\Omega}$ and all $\bt=(t_{ij})\in\mathbb{S}^3$.

Notice that the limits are functions of $\by\in\bar{\omega}$ only, that is, independent of the transversal variable $x_3$. We also recall \cite[Theorem 3.3-1]{Ciarlet4b},  which provides the asymptotic behavior of Christoffel's symbols $\Gamma_{ij}^p(\var)$, $g^{ij}(\var)$ and $g(\var)$. Indeed, if $\btheta\in\mathcal{C}^3(\bar{\omega};\mathbb{R}^3)$, then
%
\begin{align}
&\Gamma_{\alpha\beta}^\sigma(\var)=\Gamma_{\alpha\beta}^\sigma -\var x_3b_\beta^\sigma|_\alpha + {\mathcal{O} }(\var^2),\quad
\d_3 \Gamma_{\alpha\beta}^p(\var)={\mathcal{O} }(\var),\quad \Gamma_{\alpha3}^3(\var)=\Gamma_{33}^p(\var)=0,\label{eq52b1}\\
&\Gamma_{\alpha\beta}^3(\var)=b_{\alpha\beta} - \var x_3 b_\alpha^\sigma b_{\sigma\beta},\quad
\Gamma_{\alpha3}^\sigma(\var)=-b_\alpha^\sigma - \var x_3 b_\alpha^\tau b_\tau^\sigma + {\mathcal{O} }(\var^2),\\
&g^{\alpha\beta}(\var)=a^{\alpha\beta}+2\varepsilon x_3 a^{\alpha\sigma}b_\sigma^\beta+{\mathcal{O}}(\varepsilon^2),\quad g^{i3}(\varepsilon)=\delta^{i3},\quad
g(\varepsilon)=a + {\mathcal{O} }(\varepsilon),\label{eq52b2}
\end{align}
for all $\var$, $0<\var\leq\var_0$, where the order symbols ${\mathcal{O} }(\var)$ and ${\mathcal{O} }(\var^2)$  are meant with respect to the norm $\left\|\cdot\right\|_{0,\infty,\bar{\Omega}}$ defined by
\begin{align*} 
\left\|w\right\|_{0,\infty,\bar{\Omega}}=\sup \{|w(\bx)|; \bx\in\bar{\Omega}\},
\end{align*}
and the covariant derivatives $b_\beta^\sigma|_\alpha$ are defined by
$$
b_\beta^\sigma|_\alpha:=\d_\alpha b_\beta^\sigma +\Gamma^\sigma_{\alpha\tau}b_\beta^\tau - \Gamma^\tau_{\alpha\beta}b^\sigma_\tau.
$$
The functions $b_{\alpha\beta}, b_\alpha^\sigma, \Gamma_{\alpha\beta}^\sigma, b_\beta^\sigma|_\alpha$ and $a$ are identified with functions in $\mathcal{C}^0(\bar{\Omega})$. Further, there exist constants $a_0, g_0$ and $g_1$ such that
 \begin{align}
 & 0<a_0\leq a(\by) \ \forall \by\in \bar{\omega},\nonumber \\
 & 0<g_0\leq g(\varepsilon)(\bx) \leq g_1 \ \forall \bx\in\bar{\Omega} \ \textrm{and} \ \forall \ \var, 0<\varepsilon\leq \varepsilon_0.\label{g_acotado}
 \end{align}
Let the scaled heat source $q(\varepsilon):\Omega\times(0,T)\to \mathbb{R}$ and scaled applied forces $\fb(\varepsilon):\Omega\times(0,T)\to \mathbb{R}^3$ and  $\bh(\varepsilon):\Gamma_N\times(0,T)\to\mathbb{R}^3$ be defined by
\begin{align*}
&q^\var(\bx^\varepsilon)=:q(\var)(\bx)
  \quad \forall \bx\in\Omega, \ \textrm{where} \ \bx^\varepsilon=\pi^\varepsilon(\bx)\in \Omega^\varepsilon, \\
&\fb^\var=(f^{i,\varepsilon})(\bx^\varepsilon)=:\fb(\var)= (f^i(\varepsilon))(\bx)
  \quad \forall \bx\in\Omega, \ \textrm{where} \ \bx^\varepsilon=\pi^\varepsilon(\bx)\in \Omega^\varepsilon, \\
&\bh^\var=(h^{i,\varepsilon})(\bx^\varepsilon)=:\bh(\var)= (h^i(\varepsilon))(\bx)
   \quad \forall \bx\in\Gamma_N, \ \textrm{where} \ \bx^\varepsilon=\pi^\varepsilon(\bx)\in \Gamma_N^\varepsilon .
\end{align*}
Also, we define the spaces
\begin{align*} 
&V(\Omega)=\{\bv=(v_i)\in [H^1(\Omega)]^3; \bv=\mathbf{0} \ on \ \Gamma_0\},
&{S}(\Omega)=\{\varphi\in H^1(\Omega); \varphi=0 \ on \ \Gamma_0\},
\end{align*}
which are Hilbert spaces, with associated norms denoted by $\left\|\cdot\right\|_{1,\Omega}$. The scaled variational problem  can then  be written as follows:

\begin{problem}\label{problema_escalado}
Find a pair $t\mapsto(\bu(\varepsilon)(\bx,t),\vartheta(\varepsilon)(\bx,t))$ of $[0,T]\to V(\Omega)\times {S}(\Omega)$ verifying
\begin{align}
&\int_{\Omega}\rho(\ddot{u}_\alpha(\varepsilon) g^{\alpha\beta}(\var)v_\beta+\ddot{u}_3(\varepsilon) v_3)\sqrt{g(\varepsilon)}\,d{x}\nonumber%
+\int_{\Omega}A^{ijkl}(\varepsilon)e_{k||l}(\varepsilon;\bu(\varepsilon))e_{i||j}(\varepsilon;\bv)\sqrt{g(\varepsilon)}dx\\%
&\quad-\int_{\Omega}\alpha_T(3\lambda+2\mu)\vartheta(\var)(\eab(\var;\bv)g^{\alpha\beta}(\var)+\edtres(\var;\bv))\sqrt{g(\varepsilon)} dx%
\nonumber \\%
&\  =\int_{\Omega} f^{i}(\varepsilon)v_i\sqrt{g(\varepsilon)} dx
+\frac{1}{\varepsilon}\int_{\Gamma_N} h^{i}(\varepsilon)v_i\sqrt{g(\varepsilon)}  d\Gamma  \  \forall \bv\in V(\Omega),\aes,\label{ec_problema_escalado}\\
&\int_{\Omega}\beta\dot{\vartheta}(\varepsilon) \varphi\sqrt{g(\varepsilon)} dx%
+\int_{\Omega}k(\partial_\alpha \vartheta(\varepsilon) g^{\alpha\beta}(\var)\partial_\beta\varphi+\frac{1}{\var^2}\partial_3\vartheta(\varepsilon)\partial_3\varphi)\sqrt{g(\varepsilon)} dx\nonumber\\%
&\qquad+\int_{\Omega}\alpha_T(3\lambda+2\mu)\varphi(\eab(\var;\dot{\bu}(\varepsilon))g^{\alpha\beta}(\var)+\edtres(\var;\dot{\bu}(\var)))\sqrt{g(\varepsilon)} dx\nonumber\\%
&\quad=\int_{\Omega}q(\varepsilon)\varphi\,\sqrt{g(\varepsilon)} dx   \quad \forall \varphi\in{S}(\Omega),\aes,
\end{align}
with $\dot{\bu}(\varepsilon)(\cdot,0)=\bu(\varepsilon)(\cdot,0)=\bzero$ and $\vartheta(\var)(\cdot,0)=0$.
\end{problem}
%
%
%
\begin{remark}
Notice that the time-dependent version of the linearized strain tensor above is well posed when we define
$$
e_{i||j}(\varepsilon;\bu(\varepsilon))(t):=e_{i||j}(\varepsilon;\bu(\varepsilon)(t)).
$$
See for example \cite{piersanti}. Further, as commented earlier, we usually omit the explicit time dependence for the sake of a shorter notation.
\end{remark}

\begin{remark}\label{cocientes_incrementales}
The unique solvability of Problem \ref{problema_escalado} for $\var>0$ small enough is similar to Problem \ref{problema_eps} and the regularity obtained for the solutions is analogue. In particular, we find $\dot{\bu}(\var)(\cdot,t)\in V(\Omega)$ and $\dot{\vartheta}(\varepsilon)(\cdot,t)\in S(\Omega)\ \aes.$
\end{remark}
%
We now present some additional results which will be used in the next section. In \cite[Theorem 3.4-1]{Ciarlet4b}, we find the following useful result:
\begin{theorem}\label{th_int_nula}
 Let $\omega$ be a domain in $\mathbb{R}^2$ with boundary $\gamma$, let $\Omega=\omega\times (-1,1)$, and let $g\in L^p(\Omega)$, $p>1$, be a function such that
 \begin{align*}
 \intO g \d_3v dx=0, \ \textrm{for all} \ v\in \mathcal{C}^{\infty}(\bar{\Omega}) \ \textrm{with} \ v=0 \on \gamma\times[-1,1].
 \end{align*}
 Then $g=0$ a.e in $\Omega$.
\end{theorem}

%

We provide here, as a standalone theorem, a result which can be found inside the proof of \cite[Theorem 4.4-1]{Ciarlet4b}.

\begin{theorem}\label{trazatapa}
Let $X(\Omega):=\{v\in L^2(\Omega);\ \partial_3 v\in L^2(\Omega)\}$ ($\partial_3 v$ being a derivative in the sense of distributions). Then, the trace $v(\cdot,z)$ of any function $v\in X(\Omega)$ is well defined as a function in $L^2(\omega)$ for all $z\in[-1,1]$ and the trace operator defined in this fashion is continuous. In particular, there exists a constant $c_1>0$ such that
\begin{align}\nonumber
\left\|v\right\|_{L^2(\Gamma_N)}\leq c_1 \left( |v|_{0,\Omega}^2 + |\d_3v|_{0,\Omega}^2\right)^{1/2}
\end{align}
for all $\bv\in X(\Omega)$. As consequence there exists a constant  $c_2>0$ such that
\begin{align}\label{continuity}
\left\|v_3\right\|_{L^2(\Gamma_N)}\leq c_2\left( \sum_{i,j}|\eij(\varepsilon;\bv)|_{0,\Omega}^2  \right)^{1/2}\ \forall \bv\in V(\Omega).
\end{align}
\end{theorem}

\section{Formal asymptotic analysis} \label{procedure} 

In this section we briefly describe the formal procedure to identify possible two-dimensional limit problems, depending on the geometry of the middle surface, the set where the boundary conditions are given, the order of the applied forces (the procedure is described in detail in \cite{Ciarlet4b} for elastic shells in the static case).
We consider scaled applied forces and heat source of the form
$$ 
\fb(\varepsilon)(\bx)=\varepsilon^m\fb^m(\bx),\ q(\var)(\bx)=\var^m q^m(\bx) \ \forall \bx\in \Omega, \quad 
\bh(\varepsilon)(\bx)=\varepsilon^{m+1}\bh^{m+1}(\bx) \ \forall  \bx\in \Gamma_N ,
$$
where $m$ is an integer number that will show the order of the volume, heat source and surface forces, respectively.
We substitute in  (\ref{ec_problema_escalado}) to obtain the following problem:
\begin{problem}\label{problema_orden_fuerzas}
Find a pair $t\mapsto(\bu(\varepsilon)(\bx,t),\vartheta(\varepsilon)(\bx,t))$ of $[0,T]\to V(\Omega)\times {S}(\Omega)$ verifying
\begin{align}
&\int_{\Omega}\rho(\ddot{u}_\alpha(\varepsilon) g^{\alpha\beta}(\var)v_\beta+\ddot{u}_3(\varepsilon) v_3)\sqrt{g(\varepsilon)}\,d{x}\nonumber%
+\int_{\Omega}A^{ijkl}(\varepsilon)e_{k||l}(\varepsilon;\bu(\varepsilon))e_{i||j}(\varepsilon;\bv)\sqrt{g(\varepsilon)}dx\\%
&\quad-\int_{\Omega}\alpha_T(3\lambda+2\mu)\vartheta(\var)(\eab(\var;\bv)g^{\alpha\beta}(\var)+\edtres(\var;\bv))\sqrt{g(\varepsilon)} dx%
\nonumber \\%
&\quad =\int_{\Omega} \var^m f^{i,m}v_i \sqrt{g(\varepsilon)} dx
+\int_{\Gamma_N} \var^{m}h^{i,m+1}v_i\sqrt{g(\varepsilon)}  d\Gamma  \  \forall \bv\in V(\Omega),\aes,\label{ecuacion_orden_fuerzas}\\
&\int_{\Omega}\beta\dot{\vartheta}(\varepsilon) \varphi\sqrt{g(\varepsilon)} dx%
+\int_{\Omega}k(\partial_\alpha \vartheta(\varepsilon) g^{\alpha\beta}(\var)\partial_\beta\varphi+\frac{1}{\var^2}\partial_3\vartheta(\varepsilon)\partial_3\varphi)\sqrt{g(\varepsilon)} dx\nonumber\\%
&\qquad+\int_{\Omega}\alpha_T(3\lambda+2\mu)\varphi(\eab(\var;\dot{\bu}(\varepsilon))g^{\alpha\beta}(\var)+\edtres(\var;\dot{\bu}(\var)))\sqrt{g(\varepsilon)} dx\nonumber\\%
&\quad=\int_{\Omega}\var^m q^m \varphi\,\sqrt{g(\varepsilon)} dx   \quad \forall \varphi\in {S}(\Omega),\aes,\label{ecuacion_orden_calor}
\end{align}
with $\dot{\bu}(\varepsilon)(\cdot,0)=\bu(\varepsilon)(\cdot,0)=\bzero$ and $\vartheta(\var)(\cdot,0)=0$.
\end{problem}

Assume that $\btheta\in\mathcal{C}^3(\bar{\omega};\mathbb{R}^3)$ and that the scaled unknowns $\bu(\varepsilon)$, $\vartheta(\var)$ admit asymptotic expansions of the form
\begin{align}
&\bu(\varepsilon)= \bu^0 + \varepsilon \bu^1 + \varepsilon^2 \bu^2 +\ldots,\label{desarrollo_asintotico}\\
&\vartheta(\var)=\vartheta^0+\var \vartheta^1+\var^2 \vartheta^2+\ldots\nonumber
\end{align}
where $\bu^0\in V(\Omega),$ $ \bu^j\in [H^1(\Omega)]^3$, $\vartheta^0\in{S}(\Omega)$, $\vartheta^j\in H^1(\Omega)$, $j\ge1$. The assumption (\ref{desarrollo_asintotico}) implies an asymptotic expansion of the scaled linear strain as follows
\begin{align*}
\eij(\var)\equiv\eij(\varepsilon;\bu(\varepsilon))&=\frac1{\varepsilon}\eij^{-1}+ \eij^0 + \varepsilon\eij^1 + \varepsilon^2\eij^2 + \varepsilon^3\eij^3+...
\end{align*}
where,
\begin{align*}
\left\{\begin{aligned}[c]
\eab^{-1}&=0, \\
  \eatres^{-1}&=\frac{1}{2}\d_3u_\alpha^0,
   \\
  \edtres^{-1}&=\d_3u_3^0,
\end{aligned}\right.
\qquad \qquad \qquad
\left\{\begin{aligned}[c]
 \eab^0&=\frac{1}{2}(\d_\beta u_\alpha^0 + \d_\alpha u_\beta^0) - \Gamma_{\alpha\beta}^\sigma u_\sigma^0 - b_{\alpha\beta}u_3^0,
 \\
\eatres^0&=\frac{1}{2}(\d_3 u_\alpha^1 + \d_\alpha u_3^0) +   b_{\alpha}^\sigma u_\sigma^0,
\\
\edtres^0&=\d_3u_3^1,
\end{aligned}\right.\qquad
\end{align*}
\begin{align*}
\left\{\begin{aligned}[c]
\eab^1&=\frac{1}{2}(\d_\beta u_\alpha^1 + \d_\alpha u_\beta^1) - \Gamma_{\alpha\beta}^\sigma u_\sigma^1 - b_{\alpha\beta}u_3^1 + x_3(b_{\beta|\alpha}^\sigma u_\sigma^0  + b_\alpha^\sigma b_{\sigma\beta}u_3^0), \\
\eatres^1&=\frac{1}{2}(\d_3 u_\alpha^2 + \d_\alpha u_3^1) +   b_{\alpha}^\sigma u_\sigma^1 + x_3b_\alpha^\tau  b_\tau^\sigma u_\sigma^0, \\
\edtres^1&=\d_3 u_3^2.
\end{aligned}\right. \qquad \quad \qquad
\end{align*}
Besides, the functions $\eij(\varepsilon;\bv) $ admit the following expansion,
\begin{align*}
\eij(\varepsilon;\bv)=\frac{1}{\varepsilon}\eij^{-1}(\bv) + \eij^0(\bv) + \varepsilon\eij^1(\bv)+...
\end{align*}
where,
\begin{align*}
\left\{\begin{aligned}[c]
\eab^{-1}(\bv)&=0,\\
 \eatres^{-1}(\bv)&=\frac{1}{2}\d_3v_\alpha,
   \\
\edtres^{-1}(\bv)&=\d_3v_3,
\end{aligned}\right.
\qquad \qquad \quad
\left\{\begin{aligned}[c]
 \eab^0(\bv)&=\frac{1}{2}(\d_\beta v_\alpha + \d_\alpha v_\beta) - \Gamma_{\alpha\beta}^\sigma v_\sigma - b_{\alpha\beta}v_3,
 \\
\eatres^0(\bv)&=\frac{1}{2} \d_\alpha v_3 +   b_{\alpha}^\sigma v_\sigma,
\\
\edtres^0(\bv)&=0,
\end{aligned}\right.
\end{align*}
\begin{align*}
\left\{\begin{aligned}[c]
\eab^1(\bv)&=  x_3b_{\beta|\alpha}^\sigma v_\sigma  + x_3b_\alpha^\sigma b_{\sigma\beta}v_3, \\\nonumber
\eatres^1(\bv)&= x_3b_\alpha^\tau b_\tau^\sigma v_\sigma, \\\nonumber
\edtres^1(\bv)&=0.
\end{aligned}\right. \qquad \quad \qquad \qquad \qquad \qquad \qquad \qquad \qquad \qquad
\end{align*}
Upon substitution on (\ref{ecuacion_orden_fuerzas}) and (\ref{ecuacion_orden_calor}), we proceed to characterize the terms involved in the asymptotic expansions by considering different values for $m$ and grouping terms of the same order. In this way, taking in (\ref{ecuacion_orden_fuerzas}) the order $m=-2$ and particular cases of test functions, we reason that $\fb^{-2}=\bh^{-1}=\bzero$, 
which leads to $\partial_3\bu^0=0$. From (\ref{ecuacion_orden_calor}), we reason that $q^{-2}=0$ and find that $\partial_3 \vartheta^0=0$. Thus the zeroth order terms of both unknowns would be independent of the transversal variable $x_3$. Particularly, $\bu^0$ can be identified with a function $\bxi^0\in V(\omega)$, and $\vartheta^0$ can be identified with a function $\zeta^0\in{S}(\omega)$ where
$$
V(\omega):=\{\beeta=(\eta_i)\in[H^1(\omega)]^3 ; \eta_i=0 \ \textrm{on} \ \gamma_0 \},\quad {S}(\omega):=\{\varphi\in H^1(\omega); \varphi=0 \ \textrm{on} \ \gamma_0 \}.
$$
Taking $m=-1$, and using particular cases of test functions, we reason that $\fb^{-1}=\bh^0=\bzero$ 
and we find that
$$
\eatres^0=0,\quad \lambda a^{\alpha\beta}\eab^0+(\lambda+2\mu)\edtres^0=\alpha_T(3\lambda+2\mu)\vartheta^0,\quad \eab^0=\gab(\bxi^0),
$$
where
\begin{equation}\label{def_gab}
\gab(\beeta):= \frac{1}{2}(\d_\beta\eta_\alpha + \d_\alpha\eta_\beta) - \Gamma_{\alpha\beta}^\sigma\eta_\sigma -  b_{\alpha\beta}\eta_3,
\end{equation}
denote the covariant components of the linearized change of metric tensor associated with a displacement field $\eta_i\ba^i$ of the surface $S$.  From (\ref{ecuacion_orden_calor}) we reason that $q^{-1}=0$ and find that $\partial_3 \vartheta^1=0$.

Having these results in mind, for $m=0$, developing $A^{ijkl}(0)$ and taking $\bv=\beeta\in V(\omega)$ and $\varphi\in {S}(\omega)$ leads to the following two-dimensional problem, to which we shall refer as {\em thermoelastic membrane problem}:
\begin{problem}\label{problema_ab2d}
Find a pair $t\mapsto(\bxi^0(\by,t),\zeta^0(\by,t))$ of $[0,T]\to V(\omega)\times {S}(\omega)$ verifying
\begin{align*}\nonumber
&2\int_\omega\rho(\ddot{\xi}^0_\alpha a^{\alpha\beta}\eta_\beta+\ddot{\xi}^0_3\eta_3)\sqrt{a}dy%
+\int_{\omega} \a\gst(\bxi^0)\gab(\beeta)\sqrt{a}dy%
-4\int_{\omega}\frac{\alpha_T\mu(3\lambda+2\mu)}{\lambda+2\mu}\zeta^0a^{\alpha\beta}\gab(\beeta)\sqrt{a}dy\\
&\qquad 
=\int_{\omega}F^{i,0}\eta_i\sqrt{a}dy \quad \forall \beeta=(\eta_i)\in V(\omega),\aes,\\
&2\int_\omega\left(\beta+\frac{\alpha_T^2(3\lambda+2\mu)^2}{\lambda+2\mu}\right)\dot{\zeta}^0\varphi\sqrt{a}dy%
+2\int_\omega k\partial_\alpha \zeta^0a^{\alpha\beta}\partial_\beta\varphi\sqrt{a}dy\\%
&\qquad+4\int_{\omega}\frac{\alpha_T\mu(3\lambda+2\mu)}{\lambda+2\mu}\varphi a^{\alpha\beta}\gab(\dot{\bxi}^0)\sqrt{a}dy%
=\int_{\omega}Q^0\varphi\sqrt{a}dy\quad \forall \varphi\in {S}(\omega),\aes,
\end{align*}
with $\dot{\bxi}^0(\cdot,0)=\bxi^0(\cdot,0)=\bzero$ and $\zeta^0(\cdot,0)=0$.
\end{problem}
Above, we have introduced $F^{i,0}:=\int_{-1}^{1}f^{i,0}dx_3+h_N^{i,1}$, with $h_N^{i,1}(\cdot)=h^{i,1}(\cdot,+1)$, and $Q^{0}:=\int_{-1}^{1}q^0dx_3$. Also, $\a$ denotes the contravariant components of the fourth order two-dimensional elasticity tensor, defined as follows:
\begin{align} \label{tensor_a_bidimensional}
\a&:=\frac{4\lambda\mu}{\lambda + 2\mu}a^{\alpha\beta}a^{\sigma\tau} + 2\mu \ten.
\end{align}
The problem above will be analyzed in more detail in the following section. There, we shall study the existence and uniqueness of solution under additional hypotheses of geometric nature and a more suitable set of functional spaces, and provide a rigorous convergence result. To that end, the following ellipticity result for the elasticity tensor will be used. 
There exists a constant $c_e>0$ independent of the variables and $\var$,  such that
\begin{equation} \label{tensor_a_elip}
\sum_{\alpha,\beta}|t_{\alpha\beta}|^2\leq c_e \a(\by)t_{\sigma\tau}t_{\alpha\beta},
\end{equation}
for all $\by\in\bar{\omega}$ and all $\bt=(t_{\alpha\beta})\in\mathbb{S}^2$ (vector space of $2\times2$ real symmetric matrices).
\section{Elliptic membrane case. Convergence}\label{convergence}

In what follows, we assume that the family of three-dimensional linearly thermoelastic shells consist of elliptic membrane shells, that is, the middle surface of the shell $S$ is uniformly elliptic and the boundary condition of place is considered on the whole lateral face of the shell, that is, $\gamma_0=\gamma$. Further, from the formal asymptotic analysis made in the preceding section, we assume the hypotheses which led to Problem \ref{problema_ab2d}, namely,
\begin{align*}
&\fb(\varepsilon)(\bx)=\fb^0(\bx),\ q(\var)(\bx)=q^0(\bx) \ \forall \bx\in \Omega, \ 
\bh(\varepsilon)(\bx)=\varepsilon\bh^1(\bx) \ \forall  \bx\in \Gamma_N.
\end{align*}
Since there is no possible ambiguity, in what follows we drop the superindices indicating the order of the different functions.

We also recall that  for elliptic membranes it is verified the following two-dimensional Korn's type inequality (see, for example, \cite[Theorem 2.7-3]{Ciarlet4b}): there exists a constant $c_M=c_M(\omega,\btheta)>0$ such that
\begin{align} \label{Korn_elipticas}
 \left( \sum_\alpha\left\|\eta_\alpha\right\|^2_{1,\omega} + |\eta_3|_{0,\omega}^2       \right)^{1/2} \leq c_M \left( \sum_{\alpha,\beta} |\gab(\beeta)|_{0,\omega}^2  \right)^{1/2} \ \forall \beeta\in V_M(\omega),
\end{align}
where
$$
V_M(\omega):=H_0^1(\omega)\times H_0^1(\omega)\times L^2(\omega),
$$
is the right space for the well-posedness of Problem \ref{problema_ab2d}. In this section and in the sequel, $C$ represents a positive generic constant whose specific value may change from line to line, independent of $\var$ and the unknowns. Besides, for the sake of simplicity, we assume that all the parameters involved are constant. Also, the notation $\bar{\bv}$ stands for the average on $x_3$, i.e., $\bar{\bv}:=\frac{1}{2}\int_{-1}^1 \bv(x_3)dx_3$.

To favour a clearer exposition, let us reformulate Problem \ref{problema_ab2d}:
\begin{problem}\label{problema_ab}
Find a pair $t\mapsto(\bxi(\by,t),\zeta(\by,t))$ of $[0,T]\to V_M(\omega)\times H_0^1(\omega)$ verifying
\begin{align}
&2\int_\omega\rho(\ddot{\xi}_\alpha a^{\alpha\beta}\eta_\beta+\ddot{\xi}_3\eta_3)\sqrt{a}dy%
+\int_{\omega} \a\gst(\bxi)\gab(\beeta)\sqrt{a}dy%
-4\int_{\omega}\frac{\alpha_T\mu(3\lambda+2\mu)}{\lambda+2\mu}\zeta a^{\alpha\beta}\gab(\beeta)\sqrt{a}dy\nonumber\\
&\quad 
=\int_{\omega}F^{i}\eta_i\sqrt{a}dy \quad \forall \beeta=(\eta_i)\in V_M(\omega),\aes,\label{pblimite1}\\
&2\int_\omega\left(\beta+\frac{\alpha_T^2(3\lambda+2\mu)^2}{\lambda+2\mu}\right)\dot{\zeta}\varphi\sqrt{a}dy%
+2\int_\omega k\partial_\alpha \zeta a^{\alpha\beta}\partial_\beta\varphi\sqrt{a}dy\nonumber\\%
&\quad+4\int_{\omega}\frac{\alpha_T\mu(3\lambda+2\mu)}{\lambda+2\mu}\varphi a^{\alpha\beta}\gab(\dot{\bxi})\sqrt{a}dy%
=\int_{\omega}Q\varphi\sqrt{a}dy\quad \forall \varphi\in H_0^1(\omega),\aes,\label{pblimite2}
\end{align}
with $\dot{\bxi}(\cdot,0)=\bxi(\cdot,0)=\bzero$ and $\zeta(\cdot,0)=0$.
\end{problem}
Above, we have used $F^{i}:=\int_{-1}^{1}f^{i}dx_3+h_N^{i}$ with $h_N^{i}(\cdot)=h^{i}(\cdot,+1)$ and $Q:=\int_{-1}^{1}q dx_3$.
The following shows that there is a unique solution for this problem.

\begin{theorem} \label{Th_exist_unic_bid_cero}
Let $\omega$  be a domain in $\mathbb{R}^2$, let $\btheta\in\mathcal{C}^2(\bar{\omega};\mathbb{R}^3)$ be an injective mapping such that the two vectors $\ba_\alpha=\d_\alpha\btheta$ are linearly independent at all points of $\bar{\omega}$. Let $f^{i}$ and $q\in H^1(0,T;L^2(\Omega))$, $h^{i}\in H^2(0,T;L^2(\Gamma_N))$.
Then the Problem \ref{problema_ab}, has a unique solution  $(\bxi,\zeta)$ such that
\begin{align*}
&\bxi\in L^\infty(0,T;V_M(\omega)),\ \dot{\bxi}\in L^\infty(0,T;[L^2(\omega)]^3)\cap L^\infty(0,T;V_M(\omega)),\ \ddot{\bxi}\in L^\infty(0,T;[L^2(\omega)]^3),\\
&\zeta\in L^\infty(0,T;L^2(\omega))\cap L^2(0,T;H_0^1(\omega)),\ \dot{\zeta}\in L^\infty(0,T;L^2(\omega))\cap L^2(0,T;H_0^1(\omega)).
\end{align*}
\end{theorem}
Like in Theorem \ref{Th_exist_unic_3D_bid_cero}, we can cast this problem into the setting of problems solved in \cite{Francfort} or \cite[p. 359]{Marsden}, for example. Yet again, we provide an alternative proof by using the Faedo-Galerkin method.

\dem
Like in Theorem \ref{Th_exist_unic_3D_bid_cero}, we will use a Faedo-Galerkin approach to prove the existence part. Then, a proof by contradiction will show uniqueness.

{\em Existence:}
Since $V_M(\omega)$ is a separable space, there exists a countable base $\{\bz^m\}\subset V_M(\omega)$ such that
$$
V_M(\omega)=\overline{\bigcup_{m\ge1} V_m},\quad {\rm where}\ V_m=\Span\{\bz^1,\bz^2,\ldots,\bz^m\}.
$$
Similarly, there exists a countable base $\{\chi^m\}\subset H_0^1(\omega)$ such that
$$
\displaystyle H_0^1(\omega)=\overline{\bigcup_{m\ge1} S_m},\quad {\rm where}\ S_m=\Span\{\chi^1,\chi^2,\ldots,\chi^m\}.
$$
We now formulate Problem \ref{problema_ab} for the finite dimensional subspaces:
\begin{problem}\label{problema_finito}
Find a pair $t\mapsto(\bxi^m(\by,t),\zeta^m(\by,t))$ of $[0,T]\to V_m\times S_m$ verifying
\begin{align}
&2\int_\omega\rho(\ddot{\xi}^m_\alpha a^{\alpha\beta}\eta^m_\beta+\ddot{\xi}^m_3\eta^m_3)\sqrt{a}dy%
+\int_{\omega} \a\gst(\bxi^m)\gab(\beeta^m)\sqrt{a}dy%
-4\int_{\omega}\frac{\alpha_T\mu(3\lambda+2\mu)}{\lambda+2\mu}\zeta^ma^{\alpha\beta}\gab(\beeta^m)\sqrt{a}dy\nonumber\\
&\quad 
=\int_{\omega}F^{i}\eta_i^m\sqrt{a}dy \quad \forall \beeta^m=(\eta^m_i)\in V_m,\forallt,\label{pblimite1finito}\\
&2\int_\omega\left(\beta+\frac{\alpha_T^2(3\lambda+2\mu)^2}{\lambda+2\mu}\right)\dot{\zeta^m}\varphi^m\sqrt{a}dy%
+2\int_\omega k\partial_\alpha \zeta^m a^{\alpha\beta}\partial_\beta\varphi^m\sqrt{a}dy\nonumber\\%
&\quad+4\int_{\omega}\frac{\alpha_T\mu(3\lambda+2\mu)}{\lambda+2\mu}\varphi^m a^{\alpha\beta}\gab(\dot{\bxi^m})\sqrt{a}dy%
=\int_{\omega}Q\varphi^m\sqrt{a}dy\quad \forall \varphi^m\in S_m,\forallt,\label{pblimite2finito}
\end{align}
with $\dot{\bxi}^m(\cdot,0)=\bxi^m(\cdot,0)=\bzero$ and $\zeta^m(\cdot,0)=0$.
\end{problem}
Now, the classical theory of systems of ordinary differential equations guarantees the existence and uniqueness of solution for Problem \ref{problema_finito}. Taking $\beeta^m=\dot{\bxi}^m$ in (\ref{pblimite1finito}) and $\varphi^m=\zeta^m$ in (\ref{pblimite2finito}), adding both expressions and integrating the time variable in $[0,t]$ gives
\begin{align}
&\rho|\dot{\bxi}^m(t)|_{a,\omega}^2+\frac{1}{2}\|\bxi^m(t)\|_{a,\omega}^2%
+\left(\beta+\frac{\alpha_T^2(3\lambda+2\mu)^2}{\lambda+2\mu}\right)|\zeta^m(t)|_{0,\omega}^2%
+2k\int_0^t\||\zeta^m(r)\||_{a,\omega}^2\,dr\nonumber\\
&\quad
=\int_0^t\int_{\omega}Q(r)\zeta^m\sqrt{a}dy\,dr\nonumber\\
&\qquad +\int_0^t\int_{\omega}\int_{-1}^1 f^{i}(r)dx_3\dot{\xi}_i^m(r)\sqrt{a}dy\,dr%
+\int_0^t\int_{\Gamma_N}h^{i}(r)\dot{\xi}_i^m(r)\sqrt{a}d\Gamma\,dr,\label{uni_3}
\end{align}
where we have introduced the following norms:
$$
|\beeta|_{a,\omega}^2:=\int_\omega(\eta_\alpha a^{\alpha\beta}\eta_\beta+(\eta_3)^2)\sqrt{a}dy\ \forall \beeta\in [L^2(\omega)]^3,
$$
which is equivalent to the usual norm $|\cdot|_{0,\omega}$ because of the ellipticity of $(a^{\alpha\beta})$ and the regularity of $\btheta$. Also,
$$
\|\beeta\|_{a,\omega}^2:=\int_{\omega} \a\gst(\beeta)\gab(\beeta)\sqrt{a}dy\ \forall \beeta\in V_M(\omega),
$$
which is a norm in $V_M(\omega)$ because of the Korn inequality (\ref{Korn_elipticas}) and the ellipticity of $\a$ (see \eqref{tensor_a_elip}). Finally,
$$
\||\varphi\||_{a,\omega}^2:=\int_\omega \partial_\alpha \varphi a^{\alpha\beta}\partial_\beta \varphi\sqrt{a}dy,
$$
which is a norm in $H_0^1(\omega)$ equivalent to the usual $\|\cdot\|_{1,\omega}$ because of the ellipticity of $(a^{\alpha\beta})$, the regularity of $\btheta$ and the Poincar\'e inequality.

By using 
the H\"{o}lder inequality in the right-hand side terms of (\ref{uni_3}), then using Theorem \ref{trazatapa} for the terms on $\Gamma_N$ 
followed by the use of  Gronwall inequality, we obtain that the following weak convergences take place for subsequences indexed by $m$ as well:
\begin{align}
&\bxi^m\stackrel[m\to\infty]{ }{\xrightharpoonup{\ \ *\ \ }}\bxi\ {\rm in}\ L^\infty(0,T;V_M(\omega)),\quad
\dot{\bxi}^m\stackrel[m\to\infty]{ }{\xrightharpoonup{\ \ *\ \ }}\dot{\bxi}\ {\rm in}\ L^\infty(0,T;[L^2(\omega)]^3),\label{deb1}\\
&\zeta^m\stackrel[m\to\infty]{ }{\xrightharpoonup{\ \ *\ \ }}\zeta\ {\rm in}\ L^\infty(0,T;L^2(\omega)),\quad
\zeta^m\stackrel[m\to\infty]{ }{\xrightharpoonup{\ \ \ \ \ }}\zeta\ {\rm in}\ L^2(0,T;H_0^1(\omega)),\label{deb2}\\
\end{align}
Using these convergences back in (\ref{pblimite1finito})--(\ref{pblimite2finito}), 
we find 
Problem \ref{problema_ab}.

We will now prove the additional regularities for $\dot{\bxi}$, $\ddot{\bxi}$ and $\dot{\zeta}$.
First, we add equations (\ref{pblimite1finito}) and (\ref{pblimite2finito}) and write the resulting equation at times $\tilde{t}=t+h$ and $t$, with $h>0$ and $0<t\le T-h$.
Then subtract these last two equations and take $\beeta^m=\dot{\bxi}^m(\tilde{t})-\dot{\bxi}^m(t)\in V_m$
and $\varphi^m=\zeta^m(\tilde{t})-\zeta^m(t)\in S_m$
to obtain
\begin{align*}
&2\int_\omega\rho((\ddot{\xi}^m_\alpha(\tilde{t})-\ddot{\xi}^m_\alpha(t)) a^{\alpha\beta}(\dot{\xi}^m_\beta(\tilde{t})-\dot{\xi}^m_\beta(t))+ (\ddot{\xi}^m_3(\tilde{t})-\ddot{\xi}^m_3(t)) (\dot{\xi}^m_3(\tilde{t})-\dot{\xi}^m_3(t)))\sqrt{a}dy\\%
&\quad %
+\int_{\omega} \a\gst(\bxi^m(\tilde{t})-\bxi^m(t))\gab(\dot{\bxi}^m(\tilde{t})-\dot{\bxi}^m(t))\sqrt{a}dy\\%
&\quad +2\int_\omega\left(\beta+\frac{\alpha_T^2(3\lambda+2\mu)^2}{\lambda+2\mu}\right)(\dot{\zeta^m}(\tilde{t})-\dot{\zeta^m}(t)) (\zeta^m(\tilde{t})-\zeta^m(t))\sqrt{a}dy\\%
&\quad
+2\int_\omega k\partial_\alpha (\zeta^m(\tilde{t})- \zeta^m(t) )a^{\alpha\beta}\partial_\beta(\zeta^m(\tilde{t})-\zeta^m(t))\sqrt{a}dy\\%
&\quad%
=\int_{\omega}(F^{i}(\tilde{t})-F^{i}(t))(\dot{\xi}_i^m(\tilde{t})-\dot{\xi}_i^m(t))\sqrt{a}dy+ \int_{\omega}(Q(\tilde{t})-Q(t))(\zeta^m(\tilde{t})-\zeta^m(t))\sqrt{a}dy,\quad \forall t\in[0,T-h],
\end{align*}
which
 gives
\begin{align*}
&\frac{d}{dt}\rho |\dot{\bxi}^m(\tilde{t})-\dot{\bxi}^m(t)|^2_{a,\omega}
+\frac{1}{2}\frac{d}{dt}\|\bxi^m(\tilde{t})-\bxi^m(t)\|^2_{a,\omega}
+\left(\beta+\frac{\alpha_T^2(3\lambda+2\mu)^2}{\lambda+2\mu}\right)
\frac{d}{dt}|\zeta^m(\tilde{t})-\zeta^m(t)|^2_{0,\omega}\\%
&\quad
+2k\| |\zeta^m(\tilde{t})- \zeta^m(t)\| |^2_{a,\omega}
=\int_{\omega}(F^{i}(\tilde{t})-F^{i}(t))(\dot{\xi}_i^m(\tilde{t})-\dot{\xi}_i^m(t))\sqrt{a}dy\\%
&\quad%
+ \int_{\omega}(Q(\tilde{t})-Q(t))(\zeta^m(\tilde{t})-\zeta^m(t))\sqrt{a}dy,\quad \forall t\in[0,T-h].
\end{align*}
Next, we integrate in $[0,t]$ to get
\begin{align*}
&\rho |\dot{\bxi}^m(\tilde{t})-\dot{\bxi}^m(t)|^2_{a,\omega}
-\rho |\dot{\bxi}^m(h)-\dot{\bxi}^m(0)|^2_{a,\omega}
+\frac{1}{2}\|\bxi^m(\tilde{t})-\bxi^m(t)\|^2_{a,\omega}-\frac{1}{2}\|\bxi^m(h)-\bxi^m(0)\|^2_{a,\omega}\\%
&\quad
+\left(\beta+\frac{\alpha_T^2(3\lambda+2\mu)^2}{\lambda+2\mu}\right)
|\zeta^m(\tilde{t})-\zeta^m(t)|^2_{0,\omega}
-\left(\beta+\frac{\alpha_T^2(3\lambda+2\mu)^2}{\lambda+2\mu}\right)
|\zeta^m(h)-\zeta^m(0)|^2_{0,\omega}\\%
&\quad
+2k\int_{0}^{t}\| |\zeta^m(r+h)- \zeta^m(r)\| |^2_{a,\omega}dr
=\int_{0}^{t}\int_{\omega}(F^{i}(r+h)-F^{i}(r))(\dot{\xi}_i^m(r+h)-\dot{\xi}_i^m(r))\sqrt{a}\,dy\,dr\\%
&\quad%
+ \int_{0}^{t}\int_{\omega}(Q(r+h)-Q(r))(\zeta^m(r+h)-\zeta^m(r))\sqrt{a}\,dy\,dr,\quad \forall t\in[0,T-h],
\end{align*}
and dividing the equation by $h^2$ and taking limits when $h\rightarrow 0$ we obtain
\begin{align*}
&\rho |\ddot{\bxi}^m(t)|^2_{a,\omega}
+\frac{1}{2}\|\dot{\bxi}^m(t)\|^2_{a,\omega}
+\left(\beta+\frac{\alpha_T^2(3\lambda+2\mu)^2}{\lambda+2\mu}\right)
|\dot{\zeta}^m(t)|^2_{0,\omega}
+2k\int_{0}^{t}\| |\dot{\zeta}^m(r)\| |^2_{a,\omega}dr\\%
&\quad
=\rho |\ddot{\bxi}^m(0)|^2_{a,\omega}+ \frac{1}{2}\|\dot{\bxi}^m(0)\|^2_{a,\omega}+ \left(\beta+\frac{\alpha_T^2(3\lambda+2\mu)^2}{\lambda+2\mu}\right)
|\dot{\zeta}^m(0)|^2_{0,\omega} + \int_{0}^{t}\int_{\omega}\dot{F}^{i}(r)\ddot{\xi}_i^m(r)\sqrt{a}\,dy\,dr\\%
&\quad%
+ \int_{0}^{t}\int_{\omega}\dot{Q}(r)\dot{\zeta}^m(r)\sqrt{a}\,dy\,dr,\quad \forall t\in[0,T],
\end{align*}
from which, by Young's inequality, we obtain
\begin{align}\label{cotas_ace_temp_2d}
&\rho |\ddot{\bxi}^m(t)|^2_{a,\omega}
+\frac{1}{2}\|\dot{\bxi}^m(t)\|^2_{a,\omega}
+\left(\beta+\frac{\alpha_T^2(3\lambda+2\mu)^2}{\lambda+2\mu}\right)
|\dot{\zeta}^m(t)|^2_{0,\omega}
+2k\int_{0}^{t}\| |\dot{\zeta}^m(r)\| |^2_{a,\omega}dr\nonumber\\%
&\quad
\le\rho |\ddot{\bxi}^m(0)|^2_{a,\omega}+ \frac{1}{2}\|\dot{\bxi}^m(0)\|^2_{a,\omega}+ \left(\beta+\frac{\alpha_T^2(3\lambda+2\mu)^2}{\lambda+2\mu}\right)
|\dot{\zeta}^m(0)|^2_{0,\omega}\nonumber\\%
&\quad + C(\dot{\fb},\dot{\bh},\dot{q})+\bar{C}\int_{0}^{t}\left\{|\ddot{\bxi}^m(r)|_{0,\omega}^2dr
+ \int_{0}^{t}|\dot{\zeta}^m(r)|_{0,\omega}^2\right\}dr,\quad \forall t\in[0,T].
\end{align}
In order to obtain bounds for $|\ddot{\bxi}^m(0)|^2_{a,\omega}$ and $|\dot{\zeta}^m(0)|_{0,\omega}^2$ we first notice that equations (\ref{pblimite1finito}) and (\ref{pblimite2finito}) hold for $t=0$ due to the compatibility required between initial and boundary conditions. Therefore, taking $t=0$ and $\beeta^m=\ddot{\bxi}^m(0)\in V_m$ in (\ref{pblimite1finito}) and $\varphi^m=\dot{\zeta}^m(0)\in S_m$ in (\ref{pblimite2finito}) and, taking into account the initial conditions, we obtain 
\begin{align*}
&\rho|\ddot{\bxi}^m(0)|_{a,\omega}^2=\int_{\omega}F^i(0)\ddot{\bxi}^m_i(0)\,\sqrt{a}dy%
\le \frac{1}{\delta}C+\delta|\ddot{\bxi}^m(0)|_{0,\omega}^2\\
&\left(\beta+\frac{\alpha_T^2(3\lambda+2\mu)^2}{\lambda+2\mu}\right)| \dot{\zeta}^m(0)|_{0,\omega}^2 =
\int_{\omega}Q(0)\dot{\zeta}^m(0)\,\sqrt{a}\,dy   \ \le \frac{1}{\tilde{\delta}}\tilde{C}+\tilde{\delta}| \dot{\zeta}^m(0)|_{0,\omega}^2,
\end{align*}
where $\delta$ and $\tilde{\delta}$ are sufficiently small positive constants.\\
Now, back to (\ref{cotas_ace_temp_2d}), taking into account the initial conditions and the bounds above we have
\begin{align*}
&\rho |\ddot{\bxi}^m(t)|^2_{a,\omega}
+\frac{1}{2}\|\dot{\bxi}^m(t)\|^2_{a,\omega}
+\left(\beta+\frac{\alpha_T^2(3\lambda+2\mu)^2}{\lambda+2\mu}\right)
|\dot{\zeta}^m(t)|^2_{0,\omega}
+2k\int_{0}^{t}\| |\dot{\zeta}^m(r)\| |^2_{a,\omega}dr\\%
&\quad
\le C+\bar{C}\int_{0}^{t}\left\{|\ddot{\bxi}^m(r)|_{0,\omega}^2dr
+ \int_{0}^{t}|\dot{\zeta}^m(r)|_{0,\omega}^2\right\}dr,\quad \forall t\in[0,T].
\end{align*}
Next, we use the equivalence between then norms $|\cdot|_{a,\omega}$ and $|\cdot|_{0,\omega}$ and we apply Gronwall's Lemma to conclude that
\begin{equation*}
 |\ddot{\bxi}^m(t)|^2_{0,\omega}
+|\dot{\zeta}^m(t)|^2_{0,\omega}\le C,\quad \forall t\in[0,T],
\end{equation*}
and further
\begin{equation*}
\|\dot{\bxi}^m(t)\|^2_{a,\omega}
+2k\int_{0}^{t}\| |\dot{\zeta}^m(r)\| |^2_{a,\omega}dr\le C\quad \forall t\in[0,T].%
\end{equation*}
Therefore, the following weak convergences take place for subsequences still indexed by $m$.
\begin{align}
&\dot{\bxi}^m\stackrel[m\to\infty]{ }{\xrightharpoonup{\ \ *\ \ }}\dot{\bxi}\ {\rm in}\ L^\infty(0,T;V_M(\omega)),\quad
\ddot{\bxi}^m\stackrel[m\to\infty]{ }{\xrightharpoonup{\ \ *\ \ }}\ddot{\bxi}\ {\rm in}\ L^\infty(0,T;[L^2(\omega)]^3),\label{deb4}\\
&\dot{\zeta}^m\stackrel[m\to\infty]{ }{\xrightharpoonup{\ \ *\ \ }}\dot{\zeta}\ {\rm in}\ L^\infty(0,T;L^2(\omega)),\quad
\dot{\zeta}^m\stackrel[m\to\infty]{ }{\xrightharpoonup{\ \ \ \ \ }}\dot{\zeta}\ {\rm in}\ L^2(0,T;H_0^1(\omega))\label{deb5}.
\end{align}
{\em Uniqueness:} We proceed by contradiction. We first assume that there exist two solutions $(\bxi^1,\zeta^1)$ and $(\bxi^2,\zeta^2)$. Define $\bar{\bxi}=\bxi^1-\bxi^2$ and $\bar{\zeta}=\zeta^1-\zeta^2$. Now, take $\beeta=\dot{\bar{\bxi}}$ in the version of (\ref{pblimite1}) for $\bxi^1$ and $\beeta=-\dot{\bar{\bxi}}$ in the version of (\ref{pblimite1}) for $\bxi^2$.
We then sum both expresions to find that
\begin{align*}
&2\int_\omega\rho(\ddot{\bar{\xi}}_\alpha a^{\alpha\beta}\dot{\bar{\xi}}_\beta+\ddot{\bar{\xi}}_3\dot{\bar{\xi}}_3)\sqrt{a}dy%
+\int_{\omega} \a\gst(\bar{\bxi})\gab(\dot{\bar{\bxi}})\sqrt{a}dy%
-4\int_{\omega}\frac{\alpha_T\mu(3\lambda+2\mu)}{\lambda+2\mu}\bar{\zeta}a^{\alpha\beta}\gab(\dot{\bar{\bxi}})\sqrt{a}dy=0.
\end{align*}
Similarly, take $\varphi=\bar{\zeta}$ in the version of (\ref{pblimite2}) for $\zeta^1$ and $\varphi=-\bar{\zeta}$ in the version of (\ref{pblimite2}) for $\zeta^2$. Then, we sum both expresions to find that
\begin{align*}
&2\int_\omega\left(\beta+\frac{\alpha_T^2(3\lambda+2\mu)^2}{\lambda+2\mu}\right)\dot{\bar{\zeta}}\bar{\zeta}\sqrt{a}dy%
+2\int_\omega k\partial_\alpha \bar{\zeta} a^{\alpha\beta}\partial_\beta \bar{\zeta}\sqrt{a}dy%
+4\int_{\omega}\frac{\alpha_T\mu(3\lambda+2\mu)}{\lambda+2\mu}\bar{\zeta} a^{\alpha\beta}\gab(\dot{\bar{\bxi}})\sqrt{a}dy=0.
\end{align*}
Then, we add both expressions above and integrate with respect to the time variable in $[0,t]$, to find
\begin{align}
&\rho|\dot{\bar{\bxi}}(t)|_{a,\omega}^2+\frac{1}{2}\|\bar{\bxi}(t)\|_{a,\omega}^2%
+\left(\beta+\frac{\alpha_T^2(3\lambda+2\mu)^2}{\lambda+2\mu}\right)|\bar{\zeta}(t)|_{0,\omega}^2%
+2k\int_0^t\||\bar{\zeta}(r)\||_{a,\omega}^2\,dr=0.
\label{uni_1}
\end{align}
We deduce from (\ref{uni_1}) that $\bar{\bxi}=\bzero$ and $\bar{\zeta}=0$, thus showing uniqueness.

\cqd

%

Now, we present here the main result of this paper, namely that the scaled three-dimensional unknowns $(\bu(\var), \vartheta(\var))$ converge, as $\var$ tends to zero, towards a limit $(\bu,{\vartheta})$ independent of the transversal variable, and that this limit can be identified with the solution $(\bxi,\zeta)$ of the Problem \ref{problema_ab}, posed over the two-dimensional set $\omega$.

In what follows, and for the sake of simplicity, we assume that for each $\var>0$ the initial condition for the scaled linear strain is
\begin{equation} \label{condicion_inicial_def}
\eij(\var)(0,\cdot)=0,
\end{equation}
this is, the domain is on its natural state with no strains on it at the beginning of the period of observation.

\begin{theorem}\label{Th_convergencia}
Assume that $\btheta\in\mathcal{C}^3(\bar{\omega};\mathbb{R}^3)$. Consider a family of elastic elliptic shells with thickness $2\var$ approaching zero and all sharing the same elliptic middle surface $S=\btheta(\bar{\omega})$.
For all $\var$, $0<\varepsilon\leq\varepsilon_0$ let $(\bu(\varepsilon),\vartheta(\var))$ be  the solution of the associated three-dimensional scaled Problem \ref{problema_orden_fuerzas} for $m=0$. 
Then, there exist functions ${\vartheta}, u_\alpha\in H^1(\Omega)$ satisfying ${\vartheta}=0$, $u_\alpha=0$ on $\gamma\times[-1,1]$ and a function  $u_3\in L^2(\Omega)$, such that
\begin{enumerate}
\item[(a)]  $\vartheta(\varepsilon)\rightarrow {\vartheta}$, $u_\alpha(\varepsilon)\rightarrow u_\alpha$ in $H^1(\Omega)$ and $u_3(\var)\rightarrow u_3$ in $L^2(\Omega)$ when $\varepsilon\rightarrow 0$\ \aes,
\item[(b)] ${\vartheta}$ and $\bu=(u_i)$ are independent of the transversal variable $x_3$.
\end{enumerate}
Furthermore, the pair $(\bu,{\vartheta})$ can be identified with the solution of Problem \ref{problema_ab}.
\end{theorem}

\dem  We follow the structure of the proof given in \cite[Theorem 4.4-1]{Ciarlet4b} for the case of elastic elliptic membrane shells. Hence, we shall reference some steps which apply in the same manner and omit some details.
Also, for the sake of readability we may use the shorter notations $\eij(\varepsilon):=\eij(\varepsilon;\bu(\varepsilon))$. In addition to that, all references to (\ref{ecuacion_orden_fuerzas}) or (\ref{ecuacion_orden_calor}) have to be considered as for $m=0$ and drop the superindices.
The proof is divided into several parts, numbered from $(i)$ to $(v)$. 

\begin{enumerate}[\it (i)]

\item {\em \textit{A priori} boundedness and extraction of weak convergent sequences. For $\var>0$ sufficiently small, there exist bounded sequences, also indexed by $\var$, and weak limits as specified below:
\begin{align*}
&u_\alpha(\var)\stackrel[\var\to0]{ }{\xrightharpoonup{\ \ *\ \ }}u_\alpha\ {\rm in}\ L^\infty(0,T;H^1(\Omega)),\quad%
u_3(\var)\stackrel[\var\to0]{ }{\xrightharpoonup{\ \ *\ \ }}u_3\ {\rm in}\ L^\infty(0,T;L^2(\Omega)),\\%
&\dot{\bu}(\var)\stackrel[\var\to0]{ }{\xrightharpoonup{\ \ *\ \ }}\dot{\bu}\ {\rm in}\ L^\infty(0,T;[L^2(\Omega)]^3),\quad %
e_{i||j}(\varepsilon)\stackrel[\var\to0]{ }{\xrightharpoonup{\ \ *\ \ }}e_{i||j}\ {\rm in}\ L^\infty(0,T;L^2(\Omega)),\\
&\vartheta(\var)\stackrel[\var\to0]{ }{\xrightharpoonup{\ \ *\ \ }}{\vartheta}\ {\rm in}\ L^\infty(0,T;L^2(\Omega)),%
\quad\partial_\alpha \vartheta(\var)\stackrel[\var\to0]{ }{\xrightharpoonup{\ \ \ \ \ }}{\vartheta}_\alpha\ {\rm in}\ L^2(0,T;L^2(\Omega)),\\%
&\var^{-1}\partial_3 \vartheta(\var)\stackrel[\var\to0]{ }{\xrightharpoonup{\ \ \ \ \ }}{\vartheta}_{3,-1}\ {\rm in}\ L^2(0,T;L^2(\Omega)).
\end{align*}
Moreover, ${\vartheta}, u_\alpha=0$ on $\Gamma_0$.
}
\medskip

For the proof of this step we take $\bv=\dot{\bu}(\var)$ in (\ref{ecuacion_orden_fuerzas}) (see Remark \ref{cocientes_incrementales}) and $\varphi=\vartheta(\var)$ in (\ref{ecuacion_orden_calor}) and sum both expressions to find
\begin{align}
&\int_{\Omega}\rho(\ddot{u}_\alpha(\varepsilon) g^{\alpha\beta}(\var)\dot{u}_\beta(\var)+\ddot{u}_3(\varepsilon) \dot{u}_3(\var))\sqrt{g(\varepsilon)}\,d{x}%
+\int_{\Omega}A^{ijkl}(\varepsilon)e_{k||l}(\varepsilon)\dot{e}_{i||j}(\var)\sqrt{g(\varepsilon)} dx\nonumber\\%
&\qquad+\int_{\Omega}\beta\dot{\vartheta}(\varepsilon)\vartheta(\var)\sqrt{g(\varepsilon)} dx%
+\int_{\Omega}k(\partial_\alpha \vartheta(\varepsilon) g^{\alpha\beta}(\var)\partial_\beta \vartheta(\var)+\frac{1}{\var^2}\partial_3\vartheta(\varepsilon)\partial_3\vartheta(\var))\sqrt{g(\varepsilon)} dx\nonumber\\%
&\qquad =\int_{\Omega}f^{i}\dot{u}_i(\varepsilon)\sqrt{g(\varepsilon)} dx%
+\int_{\Gamma_N} h^{i}\dot{u}_i(\var)\sqrt{g(\varepsilon)}  d\Gamma
+\int_{\Omega} q \vartheta(\var)\,\sqrt{g(\varepsilon)} dx.\label{conv1}
\end{align}
We now introduce the following norms:
$$
|\bv|_{g(\var),\Omega}^2:=\int_\Omega(v_\alpha g^{\alpha\beta}(\var)v_\beta+(v_3)^2)\sqrt{g(\var)}dx\ \forall \bv\in [L^2(\Omega)]^3,
$$
which is equivalent to the usual norm $|\cdot|_{0,\Omega}$ because of the ellipticity of $(g^{\alpha\beta}(\var))$ and the regularity of $\bTheta$. Also,
$$
\|\bv\|_{A(\var),\Omega}^2:=\int_{\Omega}A^{ijkl}(\varepsilon)e_{k||l}(\varepsilon;\bv)e_{i||j}(\var;\bv)\sqrt{g(\varepsilon)} dx\ \forall \bv\in V(\Omega),
$$
which is a norm in $V(\Omega)$ because of the Korn inequality 
(see \cite[Theorem 4.4-1]{Ciarlet4b}) and the ellipticity of $A^{ijkl}(\varepsilon)$. Finally,
$$
\||\varphi\||_{g(\var),\Omega}:=\int_\Omega \partial_\alpha \varphi g^{\alpha\beta}(\var)\partial_\beta \varphi\sqrt{g(\var)}dx,
$$
which is a seminorm in ${S}(\Omega)$. Because of the uniform ellipticity of the tensors and matrices involved, and the properties of $g(\var)$, we are going to be able to use constants independent of $\var$ in the estimates below. Indeed, going back to (\ref{conv1}), we obtain
\begin{align*}
&\frac{\rho}{2}\frac{d}{dt}\{|\dot{\bu}(\var)|_{g(\var),\Omega}^2\}
+\frac{1}{2}\frac{d}{dt}\{\|\bu(\var)\|_{A(\var),\Omega}^2\}%
+\frac{\beta}{2}\frac{d}{dt}\{|\vartheta(\var)|_{0,\Omega}^2\}%
+k\||\vartheta(\var)\||_{g(\var),\Omega}^2+\frac{k}{\var^2}|\partial_3\vartheta(\var)|_{0,\Omega}^2\\%
&\quad=
\int_{\Omega}f^{i}\dot{u}_i(\varepsilon)\sqrt{g(\varepsilon)} dx%
+\int_{\Gamma_N} h^{i}\dot{u}_i(\var)\sqrt{g(\varepsilon)}  d\Gamma
+\int_{\Omega} q \vartheta(\var)\,\sqrt{g(\varepsilon)} dx.
\end{align*}
Integrating in $[0,t]$ with respect to the time variable, using the equivalences mentioned above, together with the uniformity with respect to $\var$ of the constants involved in those equivalences, integrating by parts the term with the tractions $h^i$, using Theorem \ref{trazatapa} and Young's inequality, we find that there exist a constant $C>0$ independent of $\var$ such that
\begin{align*}
&|\dot{\bu}(\var)(t)|_{0,\Omega}^2+|\eij(\var)(t)|_{0,\Omega}^2+|\vartheta(\var)(t)|_{0,\Omega}^2%
+\int_0^t(|\partial_\alpha \vartheta(\var)(r)|_{0,\Omega}^2+\frac{1}{\var^2}|\partial_3\vartheta(\var)(r)|_{0,\Omega}^2)dr\\%
&\quad
\le C(\int_0^t|\dot{\bu}(\var)(r)|_{0,\Omega}^2dr+\int_0^t|\vartheta(\var)(r)|_{0,\Omega}^2dr%
+\int_0^t|\eij(\var)(r)|_{0,\Omega}^2dr\\
&\quad+\int_0^t|\fb(r)|_{0,\Omega}^2dr+\int_0^t|q(r)|_{0,\Omega}^2dr+\int_0^t|\dot{\bh}(r)|_{0,\Gamma_N}^2dr
+|\dot{\bh}(t)|_{0,\Gamma_N}^2)
\end{align*}
Hence, by using Gronwall's inequality and the three-dimensional Korn's inequality that can be found in \cite[Theorem 4.3-1]{Ciarlet4b}, all the assertions of (i) follow.
\medskip

\item {\em The limits of the scaled unknowns, $u_i$, ${\vartheta}$ found in Step $(i)$ are independent of $x_3$.}
\medskip

The part corresponding to $u_i$ is analogous to the Step  $(ii)$ in \cite[Theorem 4.4-1]{Ciarlet4b}, so we omit it.
Regarding ${\vartheta}$, its independence on $x_3$ is a consequence of the boundedness of $\{\var^{-1}\partial_3\vartheta(\var)\}$.

\medskip
\item {\em The limits $\eij$ found in $(i)$ are independent of the variable $x_3$. Moreover, they are related with the limits $\bu:=(u_i)$ and $\vartheta$ by
\begin{align}
&\eab=\gab(\bu):= \frac{1}{2}(\d_\alpha u_\beta + \d_\beta u_\alpha) - \Gamma_{\alpha\beta}^\sigma u_\sigma -b_{\alpha\beta}u_3,\nonumber\\
&\eatres=0,\label{eatres1}\\
&\edtres= \frac{\alpha_T(3\lambda+2\mu)}{\lambda + 2\mu}{\vartheta} - \frac{\lambda}{\lambda + 2\mu} a^{\alpha \beta }\eab.\label{edtres1}
\end{align}
}
\medskip

Indeed, first considering  $\bv=\bu(\var)$ in  (\ref{eab}) and $\beeta=\bu$ in (\ref{def_gab}) (\textit{par abus de langage}, since $\bu$ is independent of $x_3$, but actually $\bu\in [H^1(\Omega)]^2\times L^2(\Omega)$),  taking into account Step  $(i)$ and the convergences $\Gamma^\sigma_{\alpha\beta}(\var)\rightarrow\Gamma^\sigma_{\alpha\beta}$ and $\Gamma^3_{\alpha\beta}(\var)\rightarrow b_{\alpha\beta} $ in $\mathcal{C}^0(\bar{\Omega})$ given by
\eqref{eq52b1}--\eqref{eq52b2}, we have that
 \begin{align}\nonumber
\eab(\varepsilon)=\frac{1}{2}(\d_\beta u_\alpha(\varepsilon) + \d_\alpha u_\beta(\varepsilon)) - \Gamma_{\alpha\beta}^p(\varepsilon)u_p(\varepsilon)\deb\eab=\gab(\bu) \ \textrm{in} \ L^2(\Omega) \aes.
    \end{align}

Moreover, $\eab$ are independent of $x_3$, as a straightforward consequence of the independence on $x_3$ of $u_i$ (Step $(ii)$). In addition, let $\bv\in V(\Omega)$. As a consequence of the definition of the scaled strains in (\ref{eab})--(\ref{edtres}), we find
\begin{align}\nonumber
 &\varepsilon\eab(\varepsilon;\bv)\rightarrow 0 \ \textrm{in} \ L^2(\Omega),\quad\varepsilon\eatres(\varepsilon;\bv)\rightarrow \frac1{2}\d_3v_\alpha \  \textrm{in} \ L^2(\Omega),
 \\\nonumber
&\varepsilon\edtres(\varepsilon;\bv)=\d_3 v_3\  \textrm{in} \ L^2(\Omega), \ \textrm{for all} \ \varepsilon>0.
\end{align}
Now, for all $\bv\in V(\Omega)$, in (\ref{ecuacion_orden_fuerzas}) we can take as test function $\var\bv\in V(\Omega)$. Then, taking into account (\ref{tensor_terminos_nulos}), we have
\begin{align*}
&\var\int_{\Omega}\rho(\ddot{u}_\alpha(\varepsilon) g^{\alpha\beta}(\var)v_\beta+\ddot{u}_3(\varepsilon) v_3)\sqrt{g(\varepsilon)}\,d{x}%
+\var\int_{\Omega}A^{ijkl}(\varepsilon)e_{k||l}(\varepsilon)e_{i||j}(\varepsilon;\bv)\sqrt{g(\varepsilon)} dx\nonumber\\
&\quad-\int_{\Omega}\alpha_T(3\lambda+2\mu)\vartheta(\var)(\var\eab(\var;\bv)g^{\alpha\beta}(\var)+\var\edtres(\var;\bv))\sqrt{g(\varepsilon)} dx%
\\
&\quad= \varepsilon\int_{\Omega} f^{i}v_i\sqrt{g(\varepsilon)} dx.%
\end{align*}
Passing to the limit as $\var\to0$, decomposing $A^{ijkl}(\var)$ into the components with different asymptotic behaviour (see \eqref{eq51b1}--\eqref{eq51b2}), the properties of $g(\varepsilon)$ (see \eqref{g_acotado}) and the convergences in Step (i), we obtain the following equality:
\begin{align}
&\int_{\Omega}\left(2\mu a^{\alpha\sigma}\eatres\d_3v_\sigma + (\lambda+2\mu)\edtres\d_3v_3\right)\sqrt{a}dx%
 +\int_{\Omega}\lambda a^{\alpha\beta}\eab\d_3v_3\sqrt{a}dx \nonumber\\
&\qquad=\int_\Omega\alpha_T(3\lambda+2\mu){\vartheta}\partial_3 v_3\sqrt{a}dx \quad \forall \bv\in V(\Omega),\aes. \label{ecuacion_integral2}
\end{align}
By taking particular test functions and using Theorem \ref{th_int_nula}, we deduce (\ref{eatres1}). Then, we go back to (\ref{ecuacion_integral2}) and use again Theorem \ref{th_int_nula} to deduce (\ref{edtres1}). The independence of $\edtres$ on $x_3$ is a consequence of this relation, as well.

\medskip

 \item {\em { We find a limit two-dimensional problem verified by functions $\bar{\bu}=(\bar{\bu}_i)$} and $\bar{\vartheta}$. In particular, since the solution of this problem is unique, the convergences on Step $(i)$ are verified for the whole families  $(\bu(\varepsilon))_{\varepsilon>0}$ and $(\vartheta(\var))_{\var>0}$. We have that  $\bar{\bu}(t)=(\bar{u}_i(t))\in V_M(\omega)$ and ${\vartheta}(t)\in {S}(\Omega)$ \aes.}

\medskip
By using \cite[Theorem 4.2-1]{Ciarlet4b}  (parts (a) and (b)), and  Step $(ii)$ we find that $\bar{u}_\alpha\in H_0^1(\omega)$ and $\bar{\vartheta}\in H_0^1(\omega)$. Therefore, $\bar{\bu}\in V_M(\omega)$ \aes. Now, let $\bv=(v_i)\in V(\Omega)$ be independent of the variable $x_3$. Then, the asymptotic behaviour of the functions  $\Gamma_{\alpha\beta}^p(\varepsilon)$ and $\Gamma_{\alpha3}^\sigma(\varepsilon)$ (see \eqref{eq52b1}--\eqref{eq52b2}) implies the following convergences when $\varepsilon\rightarrow 0$ (see (\ref{eab})--(\ref{edtres})):
\begin{align}
&\eab(\varepsilon;\bv) \rightarrow\gab(\bv):= \frac{1}{2}(\d_\alpha v_\beta + \d_\beta v_\alpha)
-\Gamma_{\alpha\beta}^\sigma v_\sigma -b_{\alpha\beta}v_3 \ \textrm{in} \ L^2(\Omega),\label{xx1}\\
&\eatres(\varepsilon;\bv)\rightarrow \frac1{2}\d_\alpha v_3 + b_\alpha^\sigma v_\sigma \ \textrm{in} \ L^2(\Omega),\quad \edtres(\varepsilon;\bv)=0.\label{xx3}
\end{align}
Having this in mind, let now $\bv=(v_i)\in V(\Omega)$ be independent of $x_3$ in (\ref{ecuacion_orden_fuerzas}) and take the limit when $\var\rightarrow 0$. In the process, we make use of the asymptotic behaviour of $A^{ijkl}(\varepsilon)$ (see \eqref{eq51b1}--\eqref{eq51b2}) and $g(\varepsilon)$ (see \eqref{g_acotado}), take into account the weak convergences $\eij(\varepsilon)\stackrel{*}{\deb}\eij$ in $L^\infty(0,T;L^2(\Omega))$,
simplify by using (\ref{eatres1}) and consider the precise limits of the functions $\eij(\varepsilon;\bv)$ in (\ref{xx1})--(\ref{xx3}). As a result, we obtain the equality
\begin{align}
&\intO\rho(\ddot{u}_\alpha a^{\alpha\beta}v_\beta+\ddot{u}_3v_3)\sqrt{a}dx%
+\intO \left( \lambda a^{\alpha \beta} a^{\sigma \tau} + \mu (a^{\alpha \sigma}a^{\beta \tau} + a^{\alpha\tau}a^{\beta\sigma})  \right) \est \gab(\bv) \sqrt{a} dx\nonumber\\ %
&\quad+\intO \lambda a^{\alpha\beta}\edtres \gab(\bv)\sqrt{a} dx%
-\intO\alpha_T(3\lambda+2\mu){\vartheta}a^{\alpha\beta}\gab(\bv)\sqrt{a}dx\nonumber\\%
&= \intO f^{i} v_i\sqrt{a} dx +\int_{\Gamma_N}h^i v_i\sqrt{a}d\Gamma\aes.
\label{ref2}
\end{align}
Using (\ref{edtres1}) and since $\bu$, $\bv$ and $\vartheta$ are all independent of $x_3$ (see Step $(ii)$), we can identify them with their averages and we obtain from (\ref{ref2}) that
\begin{align}
&2\into\rho(\ddot{\bar{u}}_\alpha a^{\alpha\beta}\bar{v}_\beta+\ddot{\bar{u}}_3\bar{v}_3)\sqrt{a}dy%
+\into \a \gst(\bar{\bu})\gab(\bar{\bv})\sqrt{a}dy%
-4\into\frac{\alpha_T\mu(3\lambda+2\mu)}{\lambda + 2\mu}{\bar{\vartheta}}a^{\alpha\beta}\gab(\bar{\bv})\sqrt{a}dy\nonumber\\%
&\quad 
=\int_{\omega}\left(\int_{-1}^{1}f^i dx_3\right)\bar{v}_i\sqrt{a} dy%
+\int_{\Gamma_N}h^i \bar{v}_i\sqrt{a}d\Gamma,\aes%
,\label{et2}
\end{align}
where  $\a$ denotes the contravariant components of the fourth order two-dimensional tensor defined in (\ref{tensor_a_bidimensional}). Now, given $\beeta=(\eta_i)\in [H_0^1(\omega)]^3$, we can define $\bv=(v_i)$ such that $\bv(\by,x_3)=\beeta(\by)$ for all $(\by,x_3)\in\Omega$. Then $\bv\in V(\Omega)$ and it is independent of $x_3$; hence, as a consequence of \cite[Theorem 4.2-1]{Ciarlet4b}, the variational problems above are satisfied for $\bar{\bv}=\beeta$.
Since both sides of the equation above are continuous linear forms with respect to  $\bar{v}_3=\eta_3\in L^2(\omega)$ for  any given $\bar{v}_\alpha\in H_0^1(\omega)$, these expressions are valid for all $\beeta=(\eta_i)\in V_M(\omega)$, since $H_0^1(\omega)$ is dense in $L^2(\omega)$.

Similarly, let $\varphi\in {S}(\Omega)$ be independent of $x_3$ in (\ref{ecuacion_orden_calor}) and take the limit when $\var\rightarrow 0$. We take into account the weak convergences in Step (i), simplify by using the time derivative of (\ref{edtres1}). As a result, we obtain the equality
\begin{align}
&2\int_\omega\left(\beta+\frac{\alpha_T^2(3\lambda+2\mu)^2}{\lambda+2\mu}\right)\dot{\bar{\vartheta}}\varphi\sqrt{a}dy%
+2\int_\omega k\partial_\alpha{\bar{\vartheta}} a^{\alpha\beta}\partial_\beta\varphi\sqrt{a}dy\nonumber\\%
&\quad+4\int_{\omega}\frac{\alpha_T\mu(3\lambda+2\mu)}{\lambda+2\mu}\varphi a^{\alpha\beta}\gab(\dot{\bar{\bu}})\sqrt{a}dy%
=\int_{\omega}Q\varphi\sqrt{a}dy\quad \forall \varphi\in H_0^1(\omega),\label{eq38b}
\end{align}
hence obtaining (\ref{pblimite2}), with $\zeta$ identified with $\bar{\vartheta}$.

\medskip

\item {\em The weak convergences are, in fact, strong.}
\medskip

For this step we first consider a case without tractions, that is, we take $\bh=\bzero$. Then we will show the changes to be made for the case with tractions. In both cases we are using a monotonicity argument. We define the quantity:
\begin{align*}
&\Lambda(\varepsilon):=%
\int_{\Omega}\rho\left((\ddot{u}_\alpha(\varepsilon)-\ddot{u}_\alpha) g^{\alpha\beta}(\var)(\dot{u}_\beta(\var)-\dot{u}_\beta)+(\ddot{u}_3(\varepsilon)-\ddot{u}_3)(\dot{u}_3(\var)-\dot{u}_3)\right)\sqrt{g(\varepsilon)}\,d{x}\nonumber\\%
&\quad+\int_{\Omega}A^{ijkl}(\varepsilon)(\ekl(\varepsilon)-\ekl)(\deij(\varepsilon)-\deij)\sqrt{g(\varepsilon)}dx\\%
&\quad+\int_{\Omega}\beta(\dot{\vartheta}(\varepsilon)-\dot{{\vartheta}})(\vartheta(\var)-{\vartheta})\sqrt{g(\varepsilon)} dx\\%
&\quad+\int_{\Omega}k\{\partial_\alpha (\vartheta(\varepsilon)-{\vartheta})g^{\alpha\beta}(\var)\partial_\beta(\vartheta(\var)-{\vartheta})%
+\frac{1}{\var^2}(\partial_3(\vartheta(\varepsilon)-{\vartheta}))^2\}\sqrt{g(\varepsilon)} dx.
\end{align*}
On one hand, we integrate with respect to the time variable in $[0,t]$ and take into account (\ref{condicion_inicial_def}) and the initial conditions in Problem \ref{problema_orden_fuerzas} to obtain
\begin{align}
&2\int_0^t\Lambda(\varepsilon)dr=%
\int_{\Omega}\rho\left((\dot{u}_\alpha(\varepsilon)-\dot{u}_\alpha) g^{\alpha\beta}(\var)(\dot{u}_\beta(\var)-\dot{u}_\beta)+(\dot{u}_3(\varepsilon)-\dot{u}_3)^2\right)\sqrt{g(\varepsilon)}\,d{x}\nonumber\\%
&\quad+\int_{\Omega}A^{ijkl}(\varepsilon)(\ekl(\varepsilon)-\ekl)({\eij}(\varepsilon)-{\eij})\sqrt{g(\varepsilon)}dx\nonumber\\%
&\quad+\int_{\Omega}\beta(\vartheta(\var)-{\vartheta})^2\sqrt{g(\varepsilon)} dx\nonumber\\%
&\quad+2\int_0^t\int_{\Omega}k\{\partial_\alpha (\vartheta(\varepsilon)-{\vartheta})g^{\alpha\beta}(\var)\partial_\beta(\vartheta(\var)-{\vartheta})%
+\frac{1}{\var^2}(\partial_3(\vartheta(\varepsilon)-{\vartheta}))^2\}\sqrt{g(\varepsilon)} dx dr,\label{intLambda}
\end{align}
and as consequence of 
(\ref{elipticidadA_eps}) and (\ref{g_acotado}), we find
\begin{align}
&\int_0^t\Lambda(\varepsilon)ds\ge C(|\dot{\bu}(\var)-\dot{\bu}|_{0,\Omega}^2+|\eij(\var)-\eij|_{0,\Omega}^2%
+|\vartheta(\var)-{\vartheta}|_{0,\Omega}^2\nonumber\\%
&\quad +\int_0^t|\partial_\alpha \vartheta(\var)-\partial_\alpha {\vartheta}|_{0,\Omega}^2 ds%
+\frac{1}{\var^2}\int_0^t|\partial_3 \vartheta(\var)-\partial_3 {\vartheta}|_{0,\Omega}^2 ds.\label{intLambda2}
\end{align}

On the other hand, from the expression of $\Lambda(\var)$ and making use of (\ref{ecuacion_orden_fuerzas})--(\ref{ecuacion_orden_calor}) for $\bv=\dot{\bu}(\varepsilon)$ and $\varphi=\vartheta(\varepsilon)$, we deduce that
\begin{align}
&\Lambda(\var)=\int_{\Omega} f^i \dot{u}_i(\varepsilon)\sqrt{g(\varepsilon)}dx%
-\frac{d}{dt}\int_{\Omega} A^{ijkl}(\varepsilon)\ekl(\varepsilon)\eij\sqrt{g(\varepsilon)}dx%
+\int_{\Omega} A^{ijkl}(\varepsilon)\ekl\deij\sqrt{g(\varepsilon)}dx\nonumber
\\%
&\quad-\frac{d}{dt}\int_{\Omega}\rho\dot{u}_\alpha(\varepsilon)g^{\alpha\beta}(\var)\dot{u}_\beta\sqrt{g(\varepsilon)}\,d{x}%
+\int_{\Omega}\rho\ddot{u}_\alpha g^{\alpha\beta}(\var)\dot{u}_\beta\sqrt{g(\varepsilon)}\,d{x}\nonumber
\\%
&\quad-\frac{d}{dt}\int_{\Omega}\rho\dot{u}_3(\varepsilon)\dot{u}_3\sqrt{g(\varepsilon)}\,d{x}%
+\int_{\Omega}\rho\ddot{u}_3\dot{u}_3\sqrt{g(\varepsilon)}\,d{x}\nonumber
\\%
&\quad+\int_{\Omega} q\vartheta(\var)\sqrt{g(\varepsilon)} dx%
-\frac{d}{dt}\int_{\Omega} \beta \vartheta(\var)\vartheta\sqrt{g(\varepsilon)} dx+\int_{\Omega}\beta \dot{\vartheta} \vartheta\sqrt{g(\varepsilon)} dx\nonumber
\\%
&\quad-\int_{\Omega}k\partial_\alpha \vartheta g^{\alpha\beta}(\var)\partial_\beta(\vartheta(\var)-\vartheta)\sqrt{g(\varepsilon)} dx%
-\int_{\Omega}k\partial_\alpha \vartheta(\var) g^{\alpha\beta}(\var)\partial_\beta \vartheta\sqrt{g(\varepsilon)} dx\nonumber
\\\label{tochote}
&\quad-\frac{1}{\var^2}\int_{\Omega}k\partial_3 \vartheta\partial_3(\vartheta(\var)-\vartheta)\sqrt{g(\varepsilon)} dx%
-\frac{1}{\var^2}\int_{\Omega}k\partial_3 \vartheta(\var)\partial_3 \vartheta\sqrt{g(\varepsilon)} dx.
\end{align}
Integrating with respect to the time variable in $[0,t]$ and taking into account the initial conditions given in Problem \ref{problema_orden_fuerzas} and (\ref{condicion_inicial_def}), we obtain
\begin{align}
&\int_0^t\Lambda(\var) dr=\int_0^t\int_{\Omega} f^i \dot{u}_i(\varepsilon)\sqrt{g(\varepsilon)}dxdr%
-\int_{\Omega} A^{ijkl}(\varepsilon)\ekl(\varepsilon)\eij\sqrt{g(\varepsilon)}dx%
+\int_0^t\int_{\Omega} A^{ijkl}(\varepsilon)\ekl\deij\sqrt{g(\varepsilon)}dxdr\nonumber
\\%
&\quad-\int_{\Omega}\rho\dot{u}_\alpha(\varepsilon)g^{\alpha\beta}(\var)\dot{u}_\beta\sqrt{g(\varepsilon)}\,d{x}%
+\int_0^t\int_{\Omega}\rho\ddot{u}_\alpha g^{\alpha\beta}(\var)\dot{u}_\beta\sqrt{g(\varepsilon)}\,d{x}dr\nonumber
\\%
&\quad-\int_{\Omega}\rho\dot{u}_3(\varepsilon)\dot{u}_3\sqrt{g(\varepsilon)}\,d{x}%
+\int_0^t\int_{\Omega}\rho\ddot{u}_3\dot{u}_3\sqrt{g(\varepsilon)}\,d{x}dr\nonumber
\\%
&\quad+\int_0^t\int_{\Omega} q\vartheta(\var)\sqrt{g(\varepsilon)} dxdr%
-\int_{\Omega} \beta \vartheta(\var)\vartheta\sqrt{g(\varepsilon)} dx+\int_0^t\int_{\Omega}\beta\dot{\vartheta} \vartheta\sqrt{g(\varepsilon)} dx dr\nonumber
\\%
&\quad-\int_0^t\int_{\Omega}k\partial_\alpha \vartheta g^{\alpha\beta}(\var)\partial_\beta(\vartheta(\var)-\vartheta)\sqrt{g(\varepsilon)} dx dr%
-\int_0^t\int_{\Omega}k\partial_\alpha \vartheta(\var) g^{\alpha\beta}(\var)\partial_\beta \vartheta\sqrt{g(\varepsilon)} dx dr\nonumber\\%
&\quad-\frac{1}{\var^2}\int_0^t\int_{\Omega}k\partial_3 \vartheta\partial_3(\vartheta(\var)-\vartheta)\sqrt{g(\varepsilon)} dxdr%
-\frac{1}{\var^2}\int_0^t\int_{\Omega}k\partial_3 \vartheta\partial_3\vartheta\sqrt{g(\varepsilon)} dx dr.\nonumber
\end{align}
Take into account that $\partial_3\vartheta=0$, and let $\varepsilon\rightarrow 0$. Using the weak convergences studied in steps $(i)$ 
and $(iv)$, the asymptotic behaviour of the functions $A^{ijkl}(\varepsilon)$ and $g(\varepsilon)$ (see \eqref{eq51b1}--\eqref{eq51b2} and \eqref{g_acotado})
and the Lebesgue dominated convergence theorem, we find that
\begin{align}
&\lim_{\varepsilon\rightarrow 0}\int_0^t\Lambda(\var)dr=\int_0^t\int_{\Omega} f^i \dot{u}_i\sqrt{a}dxdr%
-\int_0^t\intO\rho\ddot{u}_\alpha a^{\alpha\beta}\dot{u}_\beta\sqrt{a}dxdr-\int_0^t\intO\rho\ddot{u}_3\dot{u}_3\sqrt{a}dxdr%
\nonumber\\\nonumber
&\quad-\int_0^t\int_{\Omega}A^{ijkl}(0)\ekl\deij\sqrt{a}dxdr
+\int_0^t\intO q\vartheta\sqrt{a}dxdr
\\
& \quad-\int_0^t\intO\beta\dot{\vartheta}\vartheta\sqrt{a}dxdr-\int_0^t\intO k\partial_\alpha\vartheta a^{\alpha\beta}\partial_\beta \vartheta\sqrt{a}dxdr.\label{Lambda}
\end{align}
Moreover, by the expressions of $A^{ijkl}(0)$ (see \eqref{eq51b1}--\eqref{eq51b2}) and using (\ref{eatres1}) we have
\begin{align*}\nonumber
&\intO A^{ijkl}(0)\ekl\deij\sqrt{a} dx%
=\intO \left( \lambda a^{\alpha \beta} a^{\sigma \tau} + \mu (a^{\alpha \sigma}a^{\beta \tau} + a^{\alpha \tau}a^{\beta\sigma})  \right) \est \deab \sqrt{a} dx\\ \nonumber
& \qquad+ \intO \lambda a^{\alpha\beta}\edtres \deab\sqrt{a} dx
+ \intO\left(\lambda a^{\sigma\tau} \est + (\lambda + 2\mu)\edtres \right)\dedtres \sqrt{a} dx.
\end{align*}
Then, using (\ref{edtres1}), we find that (\ref{Lambda}) is actually null, since its expression above coincides with the result of adding (\ref{et2}) 
for $\bar{\bv}=\dot{\bu}$ to (\ref{eq38b}) 
for $\varphi=\vartheta$ (both integrated in $[0,t]$). Indeed,
\begin{align}
&\lim_{\varepsilon\rightarrow 0}\int_0^t\Lambda(\var)dr=\int_0^t \Big(\int_{\Omega} f^i \dot{u}_i\sqrt{a}dx 
-\intO\rho\ddot{u}_\alpha a^{\alpha\beta}\dot{u}_\beta\sqrt{a}dx-\intO\rho\ddot{u}_3\dot{u}_3\sqrt{a}dx-\frac{1}{2}\intO \a \est\deab\sqrt{a}dx\nonumber\\%
&\quad
+\intO q\vartheta\sqrt{a}dx%
-\intO\left(\beta+\frac{\alpha_T^2(3\lambda+2\mu)^2}{\lambda+2\mu}\right)\dot{\vartheta}\vartheta\sqrt{a}dx%
-\intO k\partial_\alpha\vartheta a^{\alpha\beta}\partial_\beta \vartheta\sqrt{a}dx\Big)dr=0.\label{eq42b}
\end{align}
Now, for the case where tractions are not null, in (\ref{tochote}) we have an additonal term
$$
\int_{\Gamma_N} h^i \dot{u}_i(\varepsilon)\sqrt{g(\varepsilon)}d\Gamma.
$$
We integrate (\ref{tochote}) in $[0,t]$ and integrate by parts the terms with tractions corresponding to the first two components, which can be displayed as
\begin{align}
&-\int_0^t\int_{\Gamma_N} \dot{h}^\alpha(r) u_\alpha(\varepsilon)(r)\sqrt{g(\varepsilon)}d\Gamma\,dr%
+\int_{\Gamma_N} h^\alpha(t) u_\alpha(\varepsilon)(t)\sqrt{g(\varepsilon)}d\Gamma\nonumber\\%
&\quad +\int_0^t\int_{\Gamma_N} h^3(r) \dot{u}_3(\varepsilon)(r)\sqrt{g(\varepsilon)}d\Gamma\,dr\label{tracciones}.
\end{align}
When passing to the limit $\var\to0$, the terms with $u_\alpha(\var)$ above converge by using compactness arguments, since $u_\alpha(\var)\in H^1(\Omega\times(0,T))$ and the trace into $L^2(\Gamma\times(0,T))$ is a compact operator (see \cite[p. 416]{MO}). For the term with $\dot{u}_3(\var)$ we first recall that $\dot{\bu}(\var)\in V(\Omega)$ and $\dot{\vartheta}(\var)\in S(\Omega)$ (see Remark \ref{cocientes_incrementales}). 
\medskip

Next, we use the technique of incremental coefficients in the time variable to justify additional regularity and boundedness for $\dot{\bu}(\var)$ independently of $\varepsilon$. Indeed, take the sum of both equations in Problem \ref{problema_orden_fuerzas} for $m=0$. Then consider the case for time $t+h$ and subtract the case for time $t$. Next, use $\bv=\dot{\bu}(\var)(t+h)-\dot{\bu}(\var)(t)$ and $\varphi=\vartheta(\var)(t+h)-\vartheta(\var)(t)$. We find
\begin{align*}
&\int_{\Omega}\rho((\ddot{u}_\alpha(\varepsilon)(t+h)-\ddot{u}_\alpha(\varepsilon)(t)) g^{\alpha\beta}(\var)(\dot{u}_\beta(\varepsilon)(t+h)-\dot{u}_\beta(\varepsilon)(t))\\%
&+(\ddot{u}_3(\varepsilon)(t+h)-\ddot{u}_3(\varepsilon)(t)) (\dot{u}_3(\varepsilon)(t+h)-\dot{u}_3(\varepsilon)(t)))\sqrt{g(\varepsilon)}\,d{x}\nonumber\\%
&+\int_{\Omega}A^{ijkl}(\varepsilon)(e_{k||l}(\varepsilon;\bu(\varepsilon))(t+h)-e_{k||l}(\varepsilon;\bu(\varepsilon))(t))%
(e_{i||j}(\varepsilon;\dot{\bu}(t+h))-e_{i||j}(\varepsilon;\dot{\bu}(t)))\sqrt{g(\varepsilon)}dx\\%
&+\int_{\Omega}\beta(\dot{\vartheta}(\varepsilon)(t+h)-\dot{\vartheta}(\varepsilon)(t))
({\vartheta}(\varepsilon)(t+h)-{\vartheta}(\varepsilon)(t))\sqrt{g(\varepsilon)} dx\\%
&+\int_{\Omega}k((\partial_\alpha \vartheta(\varepsilon)(t+h)-\partial_\alpha \vartheta(\varepsilon)(t)) g^{\alpha\beta}(\var)(\partial_\beta\vartheta(\var)(t+h)-\partial_\beta\vartheta(\var)(t))
+\frac{1}{\var^2}(\partial_3\vartheta(\varepsilon)(t+h)-\partial_3\vartheta(\varepsilon)(t))^2)\sqrt{g(\varepsilon)} dx\nonumber\\%
&=\int_{\Omega}(f^{i}(t+h)-f^i(t))(\dot{u}_i(\varepsilon)(t+h)-\dot{u}_i(\varepsilon)(t))\sqrt{g(\varepsilon)} dx%
+\int_{\Gamma_N}(h^{i}(t+h)-h^i(t))(\dot{u}_i(\varepsilon)(t+h)-\dot{u}_i(\varepsilon)(t))\sqrt{g(\varepsilon)}  d\Gamma\\%
&+\int_{\Omega}(q(t+h)-q(t))({\vartheta}(\varepsilon)(t+h)-{\vartheta}(\varepsilon)(t))\,\sqrt{g(\varepsilon)} dx\,\aesth.
\end{align*}
Equivalently,
\begin{align*}
&\rho\frac{1}{2}\frac{d}{dt}\left|\dot{\bu}(\varepsilon)(t+h)-\dot{\bu}(\varepsilon)(t)\right|^2_{g(\var),\Omega}%
+\frac{1}{2}\frac{d}{dt}\|e_{k||l}(\varepsilon;\bu(\varepsilon)(t+h)-\bu(\varepsilon)(t))\|^2_{{\cal A}(\var),\Omega}\\%
&+\beta\frac{1}{2}\frac{d}{dt}\left|{\vartheta}(\varepsilon)(t+h)-{\vartheta}(\varepsilon)(t)\right|^2_{0,\Omega}%
+k|\partial_\alpha \vartheta(\varepsilon)(t+h)-\partial_\alpha \vartheta(\varepsilon)(t)|_{g(\var),\Omega}^2\\%
&+\frac{k}{\var^2}|\partial_3 \vartheta(\varepsilon)(t+h)-\partial_3 \vartheta(\varepsilon)(t)|_{0,\Omega}^2\\%
&\le\int_{\Omega}(f^{i}(t+h)-f^i(t))(\dot{u}_i(\varepsilon)(t+h)-\dot{u}_i(\varepsilon)(t))\sqrt{g(\varepsilon)} dx%
+\int_{\Gamma_N}(h^{i}(t+h)-h^i(t))(\dot{u}_i(\varepsilon)(t+h)-\dot{u}_i(\varepsilon)(t))\sqrt{g(\varepsilon)}  d\Gamma\\%
&+\int_{\Omega}(q(t+h)-q(t))({\vartheta}(\varepsilon)(t+h)-{\vartheta}(\varepsilon)(t))\,\sqrt{g(\varepsilon)} dx\,\aesth.
\end{align*}
Now, integrate in $[0,t]$ and use integration by parts to find
\begin{align*}
&\rho\frac{1}{2}\left|\dot{\bu}(\varepsilon)(t+h)-\dot{\bu}(\varepsilon)(t)\right|^2_{g(\var),\Omega}%
-\rho\frac{1}{2}\left|\dot{\bu}(\varepsilon)(h)-\dot{\bu}(\varepsilon)(0)\right|^2_{g(\var),\Omega}\\%
&+\frac{1}{2}\|e_{k||l}(\varepsilon;\bu(\varepsilon)(t+h)-\bu(\varepsilon)(t))\|^2_{{\cal A}(\var),\Omega}%
-\frac{1}{2}\|e_{k||l}(\varepsilon;\bu(\varepsilon)(h)-\bu(\varepsilon)(0))\|^2_{{\cal A}(\var),\Omega}\\%
&+\beta\frac{1}{2}\left|{\vartheta}(\varepsilon)(t+h)-{\vartheta}(\varepsilon)(t)\right|^2_{0,\Omega}%
-\beta\frac{1}{2}\left|{\vartheta}(\varepsilon)(h)-{\vartheta}(\varepsilon)(0)\right|^2_{0,\Omega}\\%
&+\int_0^t\left\{k|\partial_\alpha \vartheta(\varepsilon)(r+h)-\partial_\alpha\vartheta(\varepsilon)(r)|_{g(\var),\Omega}^2%
+\frac{k}{\var^2}|\partial_3 \vartheta(\varepsilon)(r+h)-\partial_3 \vartheta(\varepsilon)(r)|_{0,\Omega}^2\right\}dr\\%
&\le\int_0^t\int_{\Omega}(f^{i}(r+h)-f^i(r))(\dot{u}_i(\varepsilon)(r+h)-\dot{u}_i(\varepsilon)(r))\sqrt{g(\varepsilon)} dxdr\\%
&+\int_0^t\int_{\Omega}(q(r+h)-q(r))({\vartheta}(\varepsilon)(r+h)-{\vartheta}(\varepsilon)(r))\,\sqrt{g(\varepsilon)} dx dr\\
&-\int_0^t\int_{\Gamma_N}(\dot{h}^{i}(r+h)-\dot{h}^i(r))({u}_i(\varepsilon)(r+h)-{u}_i(\varepsilon)(r))\sqrt{g(\varepsilon)}  d\Gamma dr\\%
&+\int_{\Gamma_N}({h}^{i}(r+h)-{h}^i(r))({u}_i(\varepsilon)(r+h)-{u}_i(\varepsilon)(r))\sqrt{g(\varepsilon)}  d\Gamma|_0^t.%
\end{align*}
Divide by $h^2$ and take limits when $h\to0$. Then,
\begin{align}
&\frac{1}{2}{\rho}\left|\ddot{{\bu}}(\varepsilon)(t)\right|^2_{g(\var),\Omega}%
+\frac{1}{2}\|e_{k||l}(\varepsilon;\dot{\bu}(\varepsilon)(t)\|^2_{{\cal A}(\var),\Omega}%
+\beta\frac{1}{2}\left|\dot{\vartheta}(\varepsilon)(t)\right|^2_{0,\Omega}\nonumber\\%
&\le\frac{1}{2}{\rho}\left|\ddot{{\bu}}(\varepsilon)(0) \right|^2_{g(\var),\Omega}%
+\frac{1}{2}\|e_{k||l}(\varepsilon;\bu(\varepsilon)(0)\|^2_{{\cal A}(\var),\Omega}%
+\beta\frac{1}{2}\left|\dot{\vartheta}(\varepsilon)(0)\right|^2_{0,\Omega}\nonumber\\%
&+\int_{0}^{t}\int_{{\Omega}}\dot{{f}}^{i}(r)\ddot{{u}}(\varepsilon)_i(r)\,\sqrt{g(\varepsilon)}  dx dr%
+\int_{0}^{t}\int_{{\Omega}}\dot{{q}}(r)\dot{{\vartheta}}(\varepsilon)(r)\,\sqrt{g(\varepsilon)} \,dx dr\nonumber\\%
&-\int_{0}^{t}\int_{{\Gamma}_N}\ddot{{h}}^{i}(r)\ddot{{u}}(\varepsilon)_i(r)\,\sqrt{g(\varepsilon)}  \,d{\Gamma} dr
+\int_{\Gamma_N}(\dot{h}^{i}(t))(\dot{u}_i(\varepsilon)(t))\sqrt{g(\varepsilon)}  d\Gamma.\label{teresa1}
\end{align}
Notice that, from the left hand side we get
$$
\|e_{k||l}(\varepsilon;\dot{\bu}(\varepsilon)(t)\|^2_{{\cal A}(\var),\Omega}%
\ge \frac{1}{2}\|e_{k||l}(\varepsilon;\dot{\bu}(\varepsilon)(t)\|^2_{{\cal A}(\var),\Omega}%
+C_1|\dot{\bu}(\varepsilon)(t)|_{0,\Omega}^2+C_2\|\dot{u}_\alpha(\var)(t)\|_{1,\Omega}^2.
$$
On the right hand side,
$$
\int_{\Gamma_N}(\dot{h}^{\alpha}(t))(\dot{u}_\alpha(\varepsilon)(t))\sqrt{g(\varepsilon)}  d\Gamma%
\le \frac{1}{\delta}\|\dot{h}^{\alpha}(t)\|_{0,\Gamma_N}^2+\delta\|\dot{u}_\alpha(\varepsilon)(t)\|_{1,\Omega}^2,
$$
and, since $\dot{\bu}(\var)\in V(\Omega)$, que can use Theorem \ref{trazatapa} to find
$$
\int_{\Gamma_N}(\dot{h}^{3}(t))(\dot{u}_3(\varepsilon)(t))\sqrt{g(\varepsilon)}  d\Gamma%
\le \frac{1}{\delta}\|\dot{h}^{3}(t)\|_{0,\Gamma_N}^2
+\delta|e_{k||l}(\varepsilon;\dot{\bu}(\varepsilon)(t)|_{0,\Omega}^2.
$$
Also,
$$
\int_{0}^{t}\int_{{\Gamma}_N}\ddot{{h}}^{\alpha}(r)\ddot{{u}}(\varepsilon)_\alpha(r)\,\sqrt{g(\varepsilon)}  \,d{\Gamma}dr\le C(\ddot{{h}})+C\int_0^t|\dot{u}_\alpha(\varepsilon)(r)|_{0,\Omega}^2dr,
$$
and 
$$
\int_{0}^{t}\int_{{\Gamma}_N}\ddot{{h}}^{3}(r)\ddot{{u}}(\varepsilon)_3(r)\,\sqrt{g(\varepsilon)}  \,d{\Gamma}dr\le C(\ddot{{h}})+C\int_0^t|e_{k||l}(\varepsilon;\dot{\bu})(\varepsilon)(t)|_{0,\Omega}^2dr,
$$
where $C(\ddot{h})$ is a constant depending on data. 

In addition to that, we can obtain estimates for $\ddot{u}(\varepsilon)(0)$ and $\dot{\vartheta}(\varepsilon)(0)$ in $L^2$ equivalent norms. Indeed, take the problem for time $t=0$. Use $\bv=\dot{\bu}(\var)(h)-\dot{\bu}(\var)(0)$ and $\varphi=\vartheta(\var)(h)-\vartheta(\var)(0)$ and the initial conditions. We also assume that $h^i(0)=0$. We find
\begin{align*}
&\int_{\Omega}\rho(\ddot{u}_\alpha(\varepsilon)(0) g^{\alpha\beta}(\var)(\dot{u}_\beta(\varepsilon)(h)-\dot{u}_\beta(\varepsilon)(0))%
+\ddot{u}_3(\varepsilon)(0)(\dot{u}_3(\varepsilon)(h)-\dot{u}_3(\varepsilon)(0)))\sqrt{g(\varepsilon)}\,d{x}\nonumber\\%
&+\int_{\Omega}\beta\dot{\vartheta}(\varepsilon)(0)
({\vartheta}(\varepsilon)(h)-{\vartheta}(\varepsilon)(0))\sqrt{g(\varepsilon)} dx\\%
&=\int_{\Omega}f^{i}(0)(\dot{u}_i(\varepsilon)(h)-\dot{u}_i(\varepsilon)(0))\sqrt{g(\varepsilon)} dx%
+\int_{\Omega}(q(0))({\vartheta}(\varepsilon)(h)-{\vartheta}(\varepsilon)(0))\,\sqrt{g(\varepsilon)} dx.
\end{align*}
Divide by $h^2$ and take limits when $h\to0$. Then,
\begin{align*}
&\frac{1}{2}{\rho}\left|\ddot{{\bu}}(\varepsilon)(0)\right|^2_{g(\var),\Omega}%
+\beta\frac{1}{2}\left|\dot{\vartheta}(\varepsilon)(0)\right|^2_{0,\Omega}\\%
&\le\int_{{\Omega}}\dot{{f}}^{i}(0)\ddot{{u}}(\varepsilon)_i(0)\,\sqrt{g(\varepsilon)}  dx%
+\int_{{\Omega}}\dot{{q}}(0)\dot{{\vartheta}}(\varepsilon)(0)\,\sqrt{g(\varepsilon)} \,dx%
\end{align*}
Now using Young's inequality we can obtain
$$
\left|\ddot{{\bu}}(\varepsilon)(0)\right|^2_{g(\var),\Omega}%
+\left|\dot{\vartheta}(\varepsilon)(0)\right|^2_{0,\Omega}\le C(f^i,q).
$$



Thus, by using Gronwall's inequality in (\ref{teresa1}), we get to obtain that $|\eij(\dot{\bu})(\var)|_{0,\Omega}^2$ is bounded independently of $\var$.

As a consequence of Step $(i)$ and (\ref{continuity}), for $v=\dot{u}_3(\var)$ we find that
$$
|\dot{u}_3(\varepsilon)|_{0,\Gamma_N}\le C |\eij(\dot{\bu})(\var)|_{0,\Omega}\ \aes.
$$
Then, there exists $\psi\in L^\infty(0,T;L^2(\Gamma_N))$ such that for a subsequence keeping the same notation, it holds  $\dot{u}_3(\varepsilon)\stackrel{*}{\deb}\psi \en L^\infty(0,T;L^2(\Gamma_N))$. Since we are in the conditions of \cite[Theorem 3.6]{ArosCao19}, 
we can identify $\psi=\dot{u}_3$.

Besides, we use Lebesgue Theorem where needed, as well. Thus, the limit of the terms with traction (\ref{tracciones}) is
\begin{align*}
&-\int_0^t\int_{\Gamma_N} \dot{h}^\alpha(r) u_\alpha(r)\sqrt{a}d\Gamma\,dr%
+\int_{\Gamma_N} h^\alpha(t) u_\alpha(t)\sqrt{a}d\Gamma\\%
&\quad+\int_0^t\int_{\Gamma_N} h^3(r) \dot{u}_3(r)\sqrt{a}d\Gamma\,dr.
\end{align*}
We can undo the integration by parts, then reason like in (\ref{eq42b}).

Finally, the  strong convergences  $\eij(\varepsilon) \rightarrow\eij $ in $L^\infty(0,T;L^2(\Omega))$  also imply the strong convergences for $u_i(\var)$, by following arguments not depending on the particular set of equations, but on arguments of differential geometry and functional analysis which do not differ from those used in \cite[Theorem 4.4-1]{Ciarlet4b}. Therefore, we just omit them and refer the interested reader to the book.


%
%
%
%
\end{enumerate}
\cqd


\section{Back to the physical framework}\label{descaling}

It remains to be proved an analogous result to the previous theorem but in terms of de-scaled unknowns. We shall present the limit problem in a de-scaled form. The scalings in Section \ref{seccion_dominio_ind} suggest the de-scalings $\xi_i^\var(\by)=\xi_i(\by)$ and $\zeta^\var(\by)=\zeta(\by)$ for all $\by\in\bar{\omega}$. This way, from Problem \ref{problema_ab} we can derive

\begin{problem}\label{problema_ab_eps}
Find a pair $t\mapsto(\bxi^\var(\by,t),\zeta^\var(\by,t))$ of $[0,T]\to V_M(\omega)\times H_0^1(\omega)$ verifying
\begin{align*}
&2\var\int_\omega\rho(\ddot{\xi}^\var_\alpha a^{\alpha\beta}\eta_\beta+\ddot{\xi}^\var_3\eta_3)\sqrt{a}dy%
+\var\int_{\omega} \aeps\gst(\bxi^\var)\gab(\beeta)\sqrt{a}dy%
\\%
&\qquad-4\var\int_{\omega}\frac{\alpha_T\mu(3\lambda+2\mu)}{\lambda+2\mu}\zeta^\var a^{\alpha\beta}\gab(\beeta)\sqrt{a}dy%
=\int_{\omega}F^{i,\var}\eta_i\sqrt{a}dy \ \forall \beeta=(\eta_i)\in V_M(\omega),\\
&2\var\int_\omega\left(\beta+\frac{\alpha_T^2(3\lambda+2\mu)^2}{\lambda+2\mu}\right)\dot{\zeta^\var}\varphi\sqrt{a}dy%
+2\var\int_\omega k\partial_\alpha^\var \zeta^\var a^{\alpha\beta}\partial_\beta^\var\varphi\sqrt{a}dy\nonumber\\%
&\qquad+4\var\int_{\omega}\frac{\alpha_T\mu(3\lambda+2\mu)}{\lambda+2\mu}\varphi a^{\alpha\beta}\gab(\dot{\bxi^\var})\sqrt{a}dy%
=\int_{\omega}Q^\var \varphi\sqrt{a}dy\quad \forall \varphi\in H_0^1(\omega),
\end{align*}
with $\dot{\bxi^\var}(\cdot,0)=\bxi^\var(\cdot,0)=\bzero$ and $\zeta^\var(\cdot,0)=0$.
\end{problem}
Above, we have used $F^{i,\var}:=\int_{-\var}^{\var}f^{i,\var}dx_3^\var +h_N^{i,\var}$, with $h_{+}^{i,\var}(\cdot)=h^{i,\var}(\cdot,\var)$, and $Q^\var=\int_{-\var}^\var q^\var dx_3^\var$. Moreover, the convergences $u_\alpha(\var) \to u_\alpha$ in $H^1(\Omega)$ and $u_3(\var) \to u_3$ in $L^2(\Omega)$ from the Theorem \ref{Th_convergencia} and \cite[Theorem 4.2-1]{Ciarlet4b} together lead to the following convergences:
$$
\frac{1}{2\var}\int^\var_{-\var} u_\alpha^\var dx_3^{\var} \to \xi_\alpha  \ \textrm{in} \ H^1(\Omega),\quad
\frac{1}{2\var}\int^\var_{-\var} u_3^\var  dx_3^{\var} \to  \xi_3 \ \textrm{in} \ L^2(\Omega),\quad
\frac{1}{2\var}\int^\var_{-\var} \zeta^\var  dx_3^{\var} \to  \zeta \ \textrm{in} \ L^2(\Omega) \aes.
$$
Furthermore, we can prove the convergences of the averages of the tangential and normal components of the three-dimensional displacement vector field. To this end, we can use the same arguments as in \cite[Theorem 4.6-1]{Ciarlet4b}.

\section{Conclusions and Outlook} \label{conclusiones}

%

We have found and mathematically justified a two-dimensional limit model for thermoelastic shells, in the particular case of the so-called elliptic membranes. To this end we used the insight provided by the asymptotic expansion method and we have justified this approach by obtaining a convergence theorem. Future work will be devoted to the asymptotic analysis of contact models, possibly thermoelastic elliptic membrane and also flexural shells, which would be found under different sets of hypotheses for the order of the functions involved or the geometry of the middle surface. We are also interested in cases when contact takes friction into account and it is coupled with other effects like  wear, adhesion or damage.

\begin{acknowledgements}
This project has received funding from the European Union's
Horizon 2020 Research and Innovation Programme under the Marie Sklodowska-Curie
Grant Agreement No 823731 CONMECH and grant MTM2016-78718-P by {\em Ministerio de Econom\'{\i}a Industria y Competitividad} of Spain with the participation of FEDER.

\end{acknowledgements}

\bibliographystyle{spmpsci}      
\bibliography{biblio_CS}   

%
%

\end{document}